\documentstyle[12pt]{article}
 \input amssym.def
\input amssym
\begin{document}
\def \Z{\Bbb Z}
\def \C{\Bbb C}
\def \R{\Bbb R}
\def \Q{\Bbb Q}
\def \N{\Bbb N}
\def \wt{{\rm wt}}
\def \tr{{\rm tr}}
\def \span{{\rm span}}
\def \Res{{\rm Res}}
\def \Res{{\rm QRes}}
\def \End{{\rm End}}
\def \E{{\rm End}}
\def \Ind {{\rm Ind}}
\def \Irr {{\rm Irr}}
\def \Aut{{\rm Aut}}
\def \Hom{{\rm Hom}}
\def \mod{{\rm mod}}
\def \ann{{\rm Ann}}
\def \<{\langle} 
\def \>{\rangle} 
\def \t{\tau }
\def \a{\alpha }
\def \e{\epsilon }
\def \l{\lambda }
\def \L{\Lambda }
\def \g{\gamma}
\def \b{\beta }
\def \om{\omega }
\def \o{\omega }
\def \c{\chi}
\def \ch{\chi}
\def \cg{\chi_g}
\def \ag{\alpha_g}
\def \ah{\alpha_h}
\def \ph{\psi_h}
\def \be{\begin{equation}\label}
\def \ee{\end{equation}}
\def \bl{\begin{lem}\label}
\def \el{\end{lem}}
\def \bt{\begin{thm}\label}
\def \et{\end{thm}}
\def \bp{\begin{prop}\label}
\def \ep{\end{prop}}
\def \br{\begin{rem}\label}
\def \er{\end{rem}}
\def \bc{\begin{coro}\label}
\def \ec{\end{coro}}
\def \bd{\begin{de}\label}
\def \ed{\end{de}}
\def \pf{{\bf Proof. }}
\def \voa{{vertex operator algebra}}

\newtheorem{thm}{Theorem}[section]
\newtheorem{prop}[thm]{Proposition}
\newtheorem{coro}[thm]{Corollary}
\newtheorem{conj}[thm]{Conjecture}
\newtheorem{exa}[thm]{Example}
\newtheorem{lem}[thm]{Lemma}
\newtheorem{rem}[thm]{Remark}
\newtheorem{de}[thm]{Definition}
\newtheorem{hy}[thm]{Hypothesis}
\makeatletter
\@addtoreset{equation}{section}
\def\theequation{\thesection.\arabic{equation}}
\makeatother
\makeatletter

\baselineskip=24pt
\begin{center}{\Large \bf An analogue of the Hom functor and 
a generalized nuclear democracy theorem}
\end{center}

\vspace{0.3cm}
\begin{center}{
Haisheng Li\\ Department of Mathematics, Rutgers University-Camden,
Camden, NJ 08102}
\end{center}

{\bf Abstract.} We give an analogue of the Hom functor and
prove a generalized form of the nuclear democracy theorem
of Tsuchiya and Kanie by using a notion of tensor product for
two modules for a vertex operator algebra.

\baselineskip=24pt

\section{Introduction}
The notion of vertex operator algebra ([B], [FHL], [FLM]) is the
algebraic counterpart of the notion of what is now usually called  ``chiral
algebra'' in conformal field theory, and
vertex operator algebra theory generalizes the theories of affine Lie algebras,
the Virasoro algebra and representations (cf. [B], [DL], [FLM], [FZ]).
It has been well known (cf. [FZ], [L1]) that the irreducible highest 
weight modules
(usually called the vacuum representations) $L(\ell,0)$ for 
an affine Lie algebra $\hat{{\bf g}}$ of level $\ell$ and $L(c,0)$ for 
the Virasoro algebra with central charge $c$ have natural vertex operator 
algebra structures. If $\ell$ is 
a positive integer, it was proved ([DL], [FZ], [L1]) that the category of 
$L(\ell,0)$-modules
 is a semi-simple category whose irreducible 
objects are irreducible highest weight integrable 
$\hat{{\bf g}}$-modules of level $\ell$ (cf. [K]).
If $c=1-\frac{6(p-q)^{2}}{pq}$, where $p,q\in \{2,3,\cdots\}$ are relatively 
prime, it was proved ([DMZ], [W]) that the  category of $L(c,0)$-modules is 
also a semi-simple category whose irreducible objects are exactly those 
irreducible Virasoro algebra modules $L(c,h)$ listed in 
[BPZ]. These give the rationality (defined in Section 2) of $L(\ell,0)$ 
and $L(c,0)$. 

To state our results, let us start with definitions 
of intertwining operator.
In the minimal models, an intertwining operator from $L(c,h_{2})$ to 
$L(c,h_{3})$ was defined in [BPZ] 
to be a primary field operator $\Phi (x)$ of weight $h$, {\it i.e.,} 
$\Phi (x)\in {\rm Hom}_{{\C}}(L(c,h_{2}),L(c,h_{3}))\{x\}$ satisfying the 
following relation:
\begin{eqnarray}\label{ep}
[L(m),\Phi(x)]=x^{m}\left( (m+1)h+x{d\over dx}\right)\Phi (x)
\end{eqnarray}
for $m\in {\Z}$. For WZW models with ${\bf g}=sl_{2}$, an intertwining operator of type 
$\left(\!\begin{array}{c}j_{3}\\j j_{2}\end{array}\!\right)$ was defined 
(cf. [TK]) as a linear map $\Phi (u,x)\in {\rm Hom}(L(\ell,j_{2}),L(\ell,j_{3}))\{x\}$ such that
(\ref{ep}) holds with $h=\frac{j(j+2)}{4(\ell+2)}$ and
\begin{eqnarray}\label{eaf}
[a(m),\Phi(u,x)]=x^{m}\Phi (au,x)\;\;\;\mbox{for }m\in {\Z}, a\in {\bf g},u\in L(j),
\end{eqnarray}
where $L(j)$ is the irreducible $sl_{2}$-module with highest weight $j$.
By employing singular vectors, Tsuchiya and Kanie proved in [TK] that 
such an intertwining operator $\Phi (\cdot,x)$ on $L(j)$
can be uniquely and naturally extended to an 
intertwining operator on $L(\ell,j)$. This is the so-called nuclear 
democracy theorem of Tsuchiya and Kanie.

On the other hand, in the context of vertex operator algebra,
an intertwining operator of type 
$\left(\!\begin{array}{c}W_{3}\\W_{1} W_{2}\end{array}\!\right)$, where 
$W_{i}$ $(i=1,2,3)$
are modules for a vertex operator algebra $V$, 
is defined
in [FHL] to be a linear map $I(\cdot,x)$ from $W_{1}$ to 
$\left({\rm Hom}(W_{2},W_{3})\right)\{x\}$ 
satisfying the $L(-1)$-bracket formula (\ref{ep}) with $m=-1$ and the 
Jacobi identity (\ref{ej}) (together with 
the truncation condition (I1) in Section 2).

An intertwining operator in the sense of
[FHL] restricted to $W_{1}(0)$ gives an intertwining operator on $W_{1}(0)$ 
in the sense of [TK] and [BPZ] for the WZW and minimal models. For WZW models,
Tsuchiya and Kanie's nuclear democracy theorem implies that the two 
definitions define the same fusion rules. The question is: do we have a 
generalized form of the nuclear
democracy theorem for an arbitrary \voa ? If $V$ is not rational, the answer is negative.
(See the appendix for a  counterexample.)
As the main result of this paper we prove a generalized form of the
nuclear democracy theorem
for a rational vertex operator algebra so that for all rational models,
the fusion rules defined in the context of vertex operator algebra coincide 
with those defined in 
the context of conformal field theory.

For WZW models, one has an affine Lie algebra $\hat{{\bf g}}$ available so 
that one can make use of the notion of Verma module and singular vectors.
To any vertex operator algebra $V$, we associate a ${\Z}$-graded Lie algebra 
$g(V)=\oplus_{n\in {\Z}}g(V)_{n}$ 
with generators $t^{n}\otimes a$ 
(linearly in $a$) for $a\in V,n\in {\Z}$ and with Borcherds'
 commutator formula (\ref{ec}) and the $L(-1)$-bracket 
formula as its defining relations (see also [B],[FFR]). 
Since $L(0)$ is a central element in $g(V)_{0}$, using the triangular 
decomposition with respect to the ${\Z}$-grading we have the notions of 
generalized Verma $g(V)$-module [Le] (or Weyl module) and lowest weight module.
Then any $V$-module $M$ is a natural $g(V)$-module such that any 
weight space $M_{(h)}$ is a natural $g(V)_{0}$-module where
 $t^{{\rm wt}a-1}\otimes a$ is represented by $a_{{\rm wt}a-1}$ for $a\in V$.
But a lowest weight, or even an irreducible lowest weight $g(V)$-module is 
not necessarily a weak $V$-module.

To formulate a nuclear democracy theorem for arbitrary rational vertex 
operator algebra, we notice that
(\ref{eaf}) is a special case of the general commutator formula (\ref{ec}).
Since (\ref{eaf}) does not hold if $a$ is not a weight-one element,
we have to use a certain cross product [FLM].
Here is our generalized form of the nuclear democracy theorem or briefly GNDT:
Let $V$ be a rational vertex operator algebra and $W_{i}$ $(i=1,2,3)$ be three
 irreducible 
$V$-modules with lowest weights $h_{i}$, respectively. Let $W_{i}(0)$ be the 
lowest weight subspace of $W_{i}$ (with weight $h_{i}$). Let $\Phi (\cdot,x)$
be a linear map from $W_{1}(0)$ to ${\rm Hom}_{{\C}}(W_{2},W_{3})\{x\}$ 
satisfying the $L(-1)$-bracket formula and
\begin{eqnarray}\label{ea}
& &(x_{1}-x_{2})^{n-1}Y(a,x_{1})\Phi(u,x_{2})
-(-x_{2}+x_{1})^{n-1}\Phi (u,x_{2})Y(a,x_{1})\nonumber\\
&=&x_{2}^{-1}\delta\left({x_{1}\over x_{2}}\right)\Phi (a_{n-1}u,x_{2})
\end{eqnarray}
for any $a\in V_{(n)},u\in W_{1}(0)$. Then there exists a unique intertwining operator 
$I(\cdot,x)$ from $W_{1}\otimes W_{2}$ to $W_{3}$ in the sense of [FHL], which extends 
$\Phi(\cdot,x)$.

To prove this GNDT, we notice that if it is true, then 
$I(\cdot,x)$ will be an injective map on $W_{1}$ so that $W_{1}(0)$ can be identified as
the space $\Phi (W_{1}(0),x)$ consisting of $\Phi (u,x)$ for $u\in U$
because $W_{1}$ is an irreducible $V$-module.
For any $u\in W_{1},a\in V$, $I(u,x)$ satisfies the $L(-1)$-bracket formula,
but
(\ref{ea}) is not true for an arbitrary $u\in W_{1}$. 
However, the local property 
holds, {i.e.,} for any $a\in V, u\in W_{1}$, there is a positive integer 
$k$ such that
$$(x_{1}-x_{2})^{k}Y(a,x_{1})\Phi (u,x_{2})=
(x_{1}-x_{2})^{k}\Phi(u,x_{2})^{k}Y(a,x_{1})$$
(cf. [DL, formula (9.37)]). 
A field operator $\Phi (x)$ from $W_{2}$ to $W_{3}$ satisfying the 
$L(-1)$-bracket formula
and the local property is called a generalized intertwining operators.
Collecting all generalized intertwining operators $\Phi(x)$ from $W_{2}$ 
to $W_{3}$ we get
a vector space $G(W_{2},W_{3})$. Then we prove (Theorem \ref{t4.6}) that 
$G(W_{2},W_{3})$ 
becomes a  $V$-module under a natural action that comes from the Jacobi 
identity. Then GNDT follows.
We also prove that $G(W_{2},W_{3})$ satisfies the universal property: For 
any $V$-module
$W$ and any intertwining operator $I(\cdot,x)$ from $W\otimes W_{2}$ to 
$W_{3}$,
there exists a unique $V$-homomorphism $\psi$ from $W$ to $G(W_{2},W_{3})$ 
such that $I(u,x)=\psi (u)(x)$ for $u\in W$. 
It follows from the universal 
property that there is a natural linear isomorphism from 
${\rm Hom}_{V}(W, G(W_{2},W_{3}))$ onto 
$I\!\left(\!\begin{array}{c}W_{3}\\W W_{2}\end{array}\!\right)$, the space
of intertwining operators of the indicated type.

For WZW models, there is another notion of intertwining 
operator involving homomorphisms from the tensor product module of
a loop module with a highest weight module to another highest weight module
for an affine Lie algebra $\hat{\frak{g}}$. By using 
the generalized form of the nuclear democracy theorem
we prove (Proposition 4.15) that 
this notion is essentially equivalent to the notion in [FHL].

The notion of $G(W_{2},W_{3})$ is clearly analogous to the notion of
``Hom''-functor. In Lie algebra theory, if $U_{i}$ $(i=1,2,3)$ are
modules for a Lie algebra ${\bf g}$, the space 
${\rm Hom}_{{\bf C}}(U_{1},U_{2})$
is a natural ${\bf g}$-module and we have the following natural 
inclusion relations:
\begin{eqnarray}
(U_{1})^{*}\otimes U_{2}\longrightarrow {\rm Hom}_{{\bf C}}(U_{1},U_{2})
\longrightarrow (U_{1}\otimes (U_{2})^{*})^{*}.
\end{eqnarray}
If both $U_{1}$ and $U_{2}$ are finite-dimensional, the arrows are 
isomorphisms so that
the space of linear homomorphisms gives a construction of tensor
 product modules.

In vertex operator algebra theory, a tensor product theory has been recently
developed [HL0-4]. (In the affine Lie algebra level, a theory of tensor 
product was
developed in [KL0-2] for modules of certain levels for an affine Lie algebra 
and part of
this theory was extended to positive integral levels in [F].)   
In [HL0-4], in addition to the notion of intertwining operator,
a notion called intertwining map was also used. 
An intertwining
map was proved to be essentially equivalent to an intertwining 
operator and could be
viewed as an operator-valued functional instead of a
formal series of operators.
As one of our results in this paper we give a definition and a construction 
of tensor product in terms of formal variable language. 

Motivated by the classical tensor product theory, we
formulate a definition of tensor product of an ordered pair of two $V$-modules
 in terms of intertwining operators and a certain universal property.
As an analogue of the construction of the classical tensor product we give 
a construction of tensor product for a rational vertex operator algebra $V$.
Roughly speaking, our tensor product module $T(W_{1},W_{2})$ is constructed 
as the 
quotient space of the tensor product vector space ${\bf 
C}[t,t^{-1}]\otimes W_{1}\otimes W_{2}$ (symbolically the linear span
of all coefficients of $Y(u_{1},x)u_{2}$ for $u_{i}\in W_{i}$)
modulo all the axioms for an intertwining operator of a certain type.
It is very natural that the tensor product vector space ${\bf
C}[t,t^{-1}]\otimes W_{1}\otimes W_{2}$ modulo all the axioms for an
intertwining operator of a certain type is a weak
$V$-module. By using universal properties, it can be proved that
the tensor product module from this 
construction is isomorphic to those (depending on $z\in {\C}^{\times}$)
constructed in [HL0-4] in the category of $V$-modules.

Analogous to the classical result, if $V$ satisfies certain ``finiteness'' and 
``semisimplicity'' conditions, we prove that there exists a unique 
maximal submodule $\Delta (W_{1},W_{2})$ inside the weak module
$G(W_{1},W_{2})$ (Proposition \ref{p4.12}) such that
$\Delta (W_{1},W_{2})'$ is a tensor product module for the ordered
pair $(W_{1},W_{2}')$ (Theorem \ref{t4.13}).

This paper is organized as follows: Section 2 is preliminary.
In Section 3 we formulate a definition of tensor product and give a 
construction of a tensor product.
In Section 4, we prove a generalized form of the nuclear democracy
theorem by using an analogue of ``Hom''- functor.   

{\bf Acknowledgment.} This paper is based on some Chapters of my
Ph.D. thesis written under 
the direction of Professors James Lepowsky and Robert Wilson at
Rutgers University, 1994. I would like to thank Professors
Lepowsky and Wilson for their insightful advice. I am greatly indebted to
Professor Lepowsky for helping me clarify many subtle points.
I am grateful to Professor Chongying Dong for
reading this paper and for giving many useful suggestions.

\newpage

\section{Preliminaries}
In this section we first review some necessary definitions from [B], [FHL] and
[FLM]. Then we present some elementary results about
certain Lie algebras and modules related to a vertex (operator) algebra. 
We use standard notations and definitions of [FHL], [FLM] and [FZ].

\bd{d2.1}
A vertex operator algebra is a quadruple $(V, Y, {\bf 1}, \omega)$ 
where $V=\oplus_{n\in {\Z}}V_{(n)}$ is a ${\Z}$-graded vector 
space, $Y(\cdot,x)$ is a linear map from
$V$ to $({\rm End} V)[[x,x^{-1}]]$, ${\bf 1}$ and $\omega$  are fixed elements
of $V$ such that the following conditions hold:

(V0)\hspace{0.25cm} $\dim V_{(n)}< \infty$ for any $n\in {\Z}$ and $V_{(n)}=0$ for $n$ 
sufficiently small;

(V1)\hspace{0.25cm} $Y({\bf 1},x)=1$;

(V2)\hspace{0.25cm} $Y(a,x){\bf 1}\in ({\rm End} V)[[x]]$ and
$\displaystyle{\lim_{x\rightarrow 0}Y(a,x){\bf 1}=a}$ for any $a\in V$;

(V3)\hspace{0.25cm} For any $a, b\in V$, $Y(a,x)b\in V((x))$ and for any 
$a, b, c\in V$,
the following {\it Jacobi identity} holds:
\begin{eqnarray}\label{ej}
& &x_{0}^{-1}\delta\left(\frac{x_{1}-x_{2}}{x_{0}}\right)Y(a,x_{1})Y(b,x_{2})c
-x_{0}^{-1}\delta\left(\frac{-x_{2}+x_{1}}{x_{0}}\right)Y(b,x_{2})Y(a,x_{1})c
\nonumber\\
&=&x_{2}^{-1}\delta\left(\frac{x_{1}-x_{0}}{x_{2}}\right)Y(Y(a,x_{0})b,x_{2})c.
\end{eqnarray}
For $a\in V$, $Y(a,x)=\sum_{n\in {\Z}}a_{n}x^{-n-1}$ is called the {\it
vertex operator} associated to $a$;

(V4) \hspace{0.25cm}Set $Y(\omega,x)=\sum_{n\in {\Z}}L(n)x^{-n-1}$. Then 
we have
\begin{eqnarray}
[L(m), L(n)]=(m-n)L(m+n)+\frac{(m^{3}-m)}{12}\delta_{m+n,0}{\rm rank }V
\end{eqnarray}
for $m,n\in {\Z}$, where ${\rm rank} V$ is a fixed complex number, called the 
{\em rank} of $V$;
\begin{eqnarray}
Y(L(-1)a,x)={d\over dx}Y(a,x)\;\;\;\mbox{ for any }a\in V;
\end{eqnarray}
and $L(0)u=nu:=({\rm wt}u)u$ for $u\in V_{(n)}, n\in {\Z}$.
\ed

This completes the definition of vertex operator algebra.

\begin{rem}\label{r2.3}
If a triple $(V, Y, {\bf 1})$ satisfies the axioms (V1)-(V3) (without assuming 
the ${\Z}$-grading and the existence of 
Virasoro algebra), $(V,Y,{\bf 1})$ is called a {\em vertex algebra}. It can 
be proved (cf. [L1]) that this
definition is equivalent to Borcherds' definition in [B]. 
\end{rem}

As a consequence of the Jacobi identity we have the following
{\em commutator formula} [B]:
\begin{eqnarray}\label{ec}
[Y(a,x_{1}),Y(b,x_{2})]={\rm
Res}_{x_{0}}x_{2}^{-1}\delta\left(\frac{x_{1}-x_{0}}{x_{2}}\right)
Y(Y(a,x_{0})b,x_{2}).
\end{eqnarray}

\bd{d2.2}
A {\it module} for a vertex operator  algebra $V$ is a pair
$(M, Y_{M})$ where $M=\oplus_{h\in {\C}}M_{(h)}$ is a ${\C}$-graded vector 
space and $Y\!_{M}\!(\cdot,x\!)\!$ is
a linear map from $V\!$ to $({\rm End} M)[[x,x^{-1}]]$ 
satisfying the following conditions:

(M0)\hspace{0.25cm} For any $h\in {\C}$, $L(0)u=hu$ for $u\in M_{(h)}$, 
$\dim M_{(h)}<\infty$ and
$M_{(n+h)}=0$ for $n\in {\Z}$ sufficiently small;

(M1)\hspace{0.25cm} $Y_{M}({\bf 1},x)=1$;
 
(M2)\hspace{0.25cm} $Y_{M}(L(-1)a,x)=\displaystyle{{d\over dx}}Y_{M}(a,x)$ 
for any $a\in V$; 
    
(M3)\hspace{0.25cm} $Y_{M}(a,x)u\in M((x))$ for any $a\in V, u\in M$ and
for any $a,b\in V, u\in M$, the following Jacobi identity holds: 
\begin{eqnarray}\label{emj}
& &x_{0}^{-1}\delta \left( \frac
{x_{1}-x_{2}}{x_{0}}\right)Y_{M}(a,x_{1})Y_{M}(b,x_{2})u-x^{-1}_{0}\delta
\left(\frac{-x_{2}+x_{1}}{x_{0}}\right)Y_{M}(b,x_{2})Y_{M}(a,x_{1})u\nonumber\\
&=&x^{-1}_{2}\delta
\left(\frac{x_{1}-x_{0}}{x_{2}}\right)Y_{M}(Y(a,x_{0})b,x_{2})u.
\end{eqnarray}
\ed

 By a {\em weak} $V$-module 
we mean a pair $(M,Y_{M})$ satisfying the axioms (M1)-(M3).
A weak $V$-module $M$ 
is said to be ${\N}$-{\em gradable} if there exists an
${\N}$-gradation $M=\oplus_{n\in{\N}}M(n)$ such that
\begin{eqnarray}
a_{n}M(k)\subseteq M(m+n-1+k)\;\;\;\mbox{ for }m,n,k\in {\Z}, a\in V_{(m)},
\end{eqnarray}
where  ${\N}$ is the set of 
nonnegative integers and $M(n)=0$ for $n<0$ by definition.
The notions of submodule, irreducible module, quotient module and module 
homomorphism can be 
defined in the obvious way. A vertex operator algebra $V$ is said to be 
{\em rational} if any ${\N}$-gradable
weak $V$-module is a direct sum of irreducible ${\N}$-gradable weak 
$V$-modules.
 If $V$ is rational,
it was proved [DLM1] that there are only finitely many irreducible modules 
up to equivalence and that
any irreducible ${\N}$-gradable weak $V$-module is a module so that
 $L(0)$ acts semisimply on any ${\N}$-gradable weak 
$V$-module. Then this definition of rationality is equivalent to Zhu's 
definition [Z] of rationality. There are also other variant 
definitions of rationality. For example, the definition of rationality
in [HL0-4] is different from the current definition.

Let $M=\oplus_{h\in {\C}}M_{(h)}$ be a $V$-module. Set 
$M'=\oplus_{h\in {\C}}M_{(h)}^{*}$ and define
\begin{eqnarray}
\< Y(a,x)u', v\>=\<u', Y(e^{xL(1)}(-x^{2})^{L(0)}a,x^{-1})v\>
\end{eqnarray}
for $u'\in M', v\in M$. Then it was proved in [FHL] that $M'$ is a $V$-module,
 called the 
{\em contragredient} module, and that $(M')'=M$.
If $f$ is a $V$-homomorphism from a $V$-module $W$ to $M$, then we have 
a $V$-homomorphism $f'$ from 
$M'$ to $W'$ such that
\begin{eqnarray}
\< f'(u'), v\>=\<u', f(v)\>\;\;\;\mbox{ for }u'\in M', v\in W.
\end{eqnarray}
Furthermore, we have $(f')'=f$ [HL0-4].

\bd{d2.5} Let $W_{1}$, $W_{2}$ and $W_{3}$ be three
weak $V$-modules. An {\it intertwining operator} of type
$\left(\!\begin{array}{c}W_{3}\\ W_{1} W_{2}\end{array}\!\right)$
is a linear map 
\begin{eqnarray}  
I(\cdot,x):& & W_{1}\rightarrow (\mbox{Hom}(W_{2},W_{3}))\{x\},\nonumber\\
& &u\mapsto I(u,x)=\sum_{\alpha \in {\C}}u_{\alpha}x^{-\alpha-1} 
\end{eqnarray}
satisfying the following conditions:

(I1)\hspace{0.25cm} For any fixed $u \in W_{1},v \in W_{2},\alpha\in
{\C}$, $u_{\alpha+n}v= 0$ for $n\in {\Z}$ sufficiently large;

(I2)\hspace{0.25cm} $I(L(-1)u,x)v=\displaystyle{{d\over dx}}I(u,x)v
\mbox{  for }u \in W_{1},v \in W_{2};$ 

(I3)\hspace{0.25cm} For $a \in V,u \in W_{1},v \in W_{2}$, the
modified Jacobi identity (\ref{ej}) where $Y(b,x_{2})$ and 
$Y(Y(a,x_{0})b,x_{2})$ are replaced
by $I(u,x_{2})$ and $I(Y(a,x_{0})u,x_{2})$, respectively, holds.
\ed

We denote by
$I\left(\!\begin{array}{c}W_{3}\\W_{1} W_{2}\end{array}\!\right)$ 
the vector space of all  
intertwining operators of the indicated type and we call the dimension of
this vector space the {\it fusion rule} of the corresponding type.


The following proposition was proved in [FHL] and [FZ]:

\begin{prop}\label{p2.6}
Let $W_{i}\!=\!\oplus_{n=0}^{\infty}W_{i}(n)$ $(i=1,2,3)$ be weak
$V\!$-modules such that $L(0)\!|_{W_{i}(n)}$ $=(h_{i}+n){\rm id}$
$(i=1,2,3)$ and let 
$I(\cdot,x)$ be an intertwining operator of type
$\left(\!\begin{array}{c}W_{3}\\ W_{1} W_{2}\end{array}\!\right)$.
Then 
\begin{eqnarray}
I^{o}(u,x):=x^{h_{1}+h_{2}-h_{3}}I(u,x)\in ({\rm
Hom}(W_{2},W_{3}))[[x,x^{-1}]]. 
\end{eqnarray}
Set $I^{o}(u,x)=\sum_{n\in {\Z}}I_{u}(n)x^{-n-1}$. Then
for any $k\in {\N}, u\in W_{1}(k),m,n \in {\N}$,
\begin{eqnarray}
I_{u}(n)W_{2}(m)\subseteq W_{3}(m+k-n-1). 
\end{eqnarray}
In particular,
\begin{eqnarray}
I_{u}(k+m+i)W_{2}(m)=0\;\;\;\mbox{ for all } i \ge 0.
\end{eqnarray}
\end{prop}

Let $W_{i}$ $(i=1,2,3)$ be $V$-modules and let $I(\cdot,x)$ be an
intertwining operator of type  
$\left(\!\begin{array}{c}W_{3}\\W_{1} W_{2}\end{array}\!\right)$. 
The {\it transpose operator} $I^{t}(\cdot,x)$ is defined by:
\begin{eqnarray}
I^{t}(\cdot,x):& &W_{2}\otimes W_{1}\rightarrow W_{3}\{x\}\nonumber\\
& &I^{t}(u_{2},x)u_{1}=e^{xL(-1)}I(u_{1},e^{\pi i}x)u_{2}
\end{eqnarray}
for $u_{1}\in W_{1},u_{2}\in W_{2}$.
The {\it adjoint operator}
$I'(\cdot,x)$ is defined by: 
\begin{eqnarray}
I'(\cdot,x)&:&W_{1}\otimes W_{3}' \rightarrow W_{2}'\{x\}\nonumber\\
& &\langle I'(u_{1},x)u_{3}',u_{2}\rangle =\langle
u_{3}',I(e^{xL(1)}(e^{\pi i}x^{-2}) 
^{L(0)}u_{1},x^{-1})u_{2}\rangle \end{eqnarray}
for $u_{1}\in W_{1},u_{2}\in W_{2},u_{3}'\in W_{3}'$. 
The following proposition was proved in [HL0-4] (see also [FHL], [L2]). 

\begin{prop}\label{p2.7}
The transpose operator $I^{t}(\cdot,x)$
and the adjoint operator $I'(\cdot,x)$ are intertwining
operators of corresponding types.
\end{prop}

Notice that the transpose operator $I^{t}(\cdot,x)$ can be defined more
generally for weak $V$-modules $W_{i}$ for $i=1,2,3$
and it follows from the same proof that it is an intertwining operator.

The following Borcherds' examples of vertex 
algebras [B] show that the notion of vertex 
algebra is really a generalization of the notion of commutative associative 
algebra.

\begin{exa}\label{e2.4}  
Let $A$ be a commutative associative algebra with identity together with a 
derivation $d$. 
Define
\begin{eqnarray}
Y(a,x)b=\left(e^{xd}a\right)b\;\;\;\;\mbox{for any }a,b\in A.
\end{eqnarray}
Then $(A,Y,1)$ is a vertex algebra.
Furthermore, let $M$ be a module for $A$ viewed as an associative algebra.
Define $Y_{M}(a,x)u=\left(e^{xd}a\right)u$ for
$a\in V,u\in M$. Then $(M,Y_{M})$ is a module for the vertex algebra
$(A,Y,1)$. In particular, let $A={\C}((t))$ and $\displaystyle{d={d\over dt}}$. Then
 $({\C}((t)),Y,1)$ is a vertex algebra.
By definition, we have
\begin{eqnarray}
Y(f(t),x)=e^{x{d\over dt}}f(t)=f(t+x)\;\;\;\mbox{ for }f(t)\in {\C}((t)).
\end{eqnarray}
It is clear that the Laurent polynomial ring  ${\C}[t,t^{-1}]$ is a
vertex subalgebra. 
\end{exa}

For convenience in the following we shall associate Lie algebras $g_{0}(V)$ and
$g(V)$ to a vertex algebra $V$. The following lemma could be found in [B]:

\begin{lem}\label{l2.8}
Let $(V,Y,{\bf 1})$ be a vertex
algebra and let $d$ be the endomorphism of $V$ defined by $d(a)=a_{-2}{\bf 1}$ 
for $a\in V$. Then the quotient space $g_{0}(V):=V/dV$ is a Lie
algebra with the bilinear product:
$[\bar{a},\bar{b}]=\overline{a_{0}b}$ for $a,b\in 
V$. Furthermore, any $V$-module $M$ is a $g_{0}(V)$-module with
the action given by: $au=a_{0}u$ for $a\in V,u\in M$.
\end{lem}

Let $V$ be any vertex algebra. Then from [FHL] (see also [B]), $\hat{V}:={\bf
C}[t,t^{-1}]\otimes V$ has a vertex algebra 
structure with $Y(f(t)\otimes u,x)=Y(f(t),x)\otimes Y(u,x)$
for any $f(t)\in {\C}[t,t^{-1}],u\in V$, and ${\bf 1}= 1\otimes {\bf
1}_{V}$. 
(The affinization of a vertex operator algebra has also been used in [HL0-4].)
 Set $\hat{d}:={d\over dt}\otimes 1+1\otimes d_{V}$. Then
$\hat{d}(u)=u_{-2}1$ for $u\in \hat{V}$.
Then from Lemma \ref{l2.8}
 $g_{0}(\hat{V})=\hat{V}/\hat{d}\hat{V}$ is a Lie algebra. For any
$m,n\in {\Z}, a \in V$, by definition we have
\begin{eqnarray}\label{e2.21}
(t^{m}\otimes a)_{n}&=&{\rm Res}_{x}x^{n}Y(t^{m}\otimes a,x)\nonumber\\
&=&{\rm Res}_{x}x^{n}(t+x)^{n}\otimes Y(a,x)\nonumber\\
&=&\sum_{i=0}^{\infty}{m\choose i}t^{m+n-i}\otimes a_{i}.
\end{eqnarray}
Thus
\begin{eqnarray}
[\overline{(t^{m}\otimes a)}, \overline{(t^{n}\otimes b)}]
=\overline{(t^{m}\otimes a)_{0}(t^{n}\otimes b)}
=\sum_{i=0}^{\infty}{m\choose i}\overline{t^{m+n-i}\otimes a_{i}b}
\end{eqnarray}
for any $a,b\in V, m,n\in {\Z}$, where ``bar'' denotes the
natural quotient map from $\hat{V}$ to $g_{0}(\hat{V})$. Therefore, we have
(see also [B])

\begin{prop}\label{p2.9}
Let $V$ be any vertex algebra. Then the
quotient space $g(V):=g_{0}(\hat{V})$ is a Lie algebra with
the bilinear operation:
\begin{eqnarray}
[\overline{t^{m}\otimes a},\overline{t^{n}\otimes b}]
=\sum_{i=0}^{\infty}{m\choose i}\overline{t^{m+n-i}\otimes a_{i}b}.
\end{eqnarray}
\end{prop}

(This Lie algebra $g(V)$ has been also studied in [FFR].) We also use $a(m)$ 
for $t^{m}\otimes a$ through the paper.
It is clear that ${\bf 1}(-1)$ is a central
element of $g(V)$. If ${\bf 1}(-1)$ acts as a scalar $k$ on a
$g(V)$-module $M$, we call $M$ a $g(V)$-module of {\it level $k$}.
(This corresponds to level for affine Lie algebras.) 
A $g(V)$-module $M$ is said to be {\em restricted} if for any $a\in V, 
u\in M$, $a(n)u=0$ for $n$ sufficiently large. Then any weak $V$-module $M$ 
is a restricted 
$g(V)$-module of level one, where $a(n)$ is represented by $a_{n}$. (However, a restricted $g(V)$-module is not necessarily a weak $V$-module.)
Then we obtain a functor $\cal{F}$ from the category of weak $V$-modules to the 
category of restricted $g(V)$-modules.
For any restricted $g(V)$-module $M$, we define $J(M)$ to be the intersection 
of all $\ker f$, where $f$ runs through all $g(V)$-homomorphisms from $M$ to weak $V$-modules. 
Then $M$ is a weak $V$-module if and only if $J(M)=0$. Furthermore, $M/J(M)$ is a weak $V$-module and $M/J(M)$ is a universal from $M$ to the functor 
$\cal{F}$ [J].

To summarize, for any vertex algebra $V$ we have two Lie algebras
$g_{0}(V)$ and $g(V)$ which are related by the following inclusion
relations: 
\begin{eqnarray}
g_{0}(V)\subseteq g(V)\simeq g_{0}(\hat{V})\subseteq g(\hat{V})\subseteq\cdots.
\end{eqnarray}

Let $V$ be a vertex operator algebra. For any $a\in
V_{(m)}, m, n\in {\Z}$, we define 
\begin{eqnarray}
{\rm deg}\:a(n)=\deg\: (t^{n}\otimes a)={\rm wt a}-n-1=m-n-1.
\end{eqnarray}
Then $g(V)$ becomes
a ${\Z}$-graded Lie algebra. Denote by $g(V)_{0}$ the degree-zero subalgebra. 
Then we obtain a triangular decomposition $g(V)=g(V)_{+}\oplus g(V)_{0}\oplus g(V)_{-}$.

\begin{lem}\label{l2.11}
Let $V$ be a vertex algebra, let $M$ be a $V$-module
and let $z$ be any nonzero complex number. For any $a\in V, u\in M,
m,n\in {\Z}$, define
\begin{eqnarray}
a(m)(t^{n}\otimes u)=\sum_{i=0}^{\infty}{m\choose i}z^{m-i}(t^{m+n-i}\otimes a_{i}u).
\end{eqnarray}
Then this defines a $g(V)$-module (of level zero) structure on
$\hat{M}:={\C}[t,t^{-1}]\otimes M$.
\end{lem}

{\bf Proof.} Let $\psi$ be the automorphism of the associative algebra
${\C}[t,t^{-1}]$ such that $\psi(f(t))=f(zt)$ for $f(t)\in {\C}[t,t^{-1}]$. 
Set
${\C}[t,t^{-1}]^{\psi}={\C}[t,t^{-1}]$. Then ${\C}[t,t^{-1}]^{\psi}$ is a ${\C}[t,t^{-1}]$-module
with the following action:
$$ f(t)u=\psi (f(t))u=f(zt)u\;\;\;\mbox{ for }f(t)\in {\C}[t,t^{-1}],u\in {\C}[t,t^{-1}]^{\psi}.$$
By Example \ref{e2.4} ${\C}[t,t^{-1}]^{\psi}$
is a module for the vertex algebra ${\C}[t,t^{-1}]$ such that
\begin{eqnarray}
Y(f(t),x)u=\psi\left(e^{x{d\over dt}}f(t)\right)u=f(zt+x)u
\end{eqnarray}
for $f(t)\in {\C}[t,t^{-1}], u\in {\C}[t,t^{-1}]^{\psi}$. 
Then ${\C}[t,t^{-1}]^{\psi}\otimes M$
is a $\hat{V}$-module, so that it is a $g(V)$ ($=g_{0}(\hat{V})$)-module
(of level zero). Then the lemma follows (\ref{e2.21}) immediately.$\;\;\;\;\Box$

Let $V$ be a vertex algebra and let $M$ be a $V$-module.
For any nonzero complex number $z$, let ${\C}_{z}$ be the
evaluation module for the associative algebra ${\C}[t,t^{-1}]$ with
$t$ acting as a scalar $z$. 
Then from Example \ref{e2.4} ${\C}_{z}$ is a module for vertex algebra
${\C}[t,t^{-1}]$, 
so that ${\C}_{z}\otimes M$ is a $\hat{V}$-module. 
Therefore ${\C}_{z}\otimes M$
is a $g(V)=g_{0}(\hat{V})$-module (by Lemma 2.1). From (\ref{e2.21}) we have
\begin{eqnarray}
a(m)\cdot (1\otimes u)=\sum_{i=0}^{\infty}{m\choose i}
z^{m-i}(1\otimes a_{i}u)\;\;\;\mbox{for }a\in V,u\in M.
\end{eqnarray}
Denote this 
$g(V)$-module by $M_{z}$. Then we obtain

\bp{p2.10}
Let $V$ be a vertex algebra, let $M$ be a $V$-module and let $z$ be any nonzero complex number.
Define $\rho: g(V)\rightarrow {\rm End}_{{\C}}M$ as follows:
\begin{eqnarray}
\rho (a(m))u=\sum_{i=0}^{\infty}{m\choose i}
z^{m-i}a_{i}u\;\;\;\mbox{for }a\in V,u\in M.
\end{eqnarray}
Then $\rho$ is a representation of $g(V)$ (of level zero) on $M$. 
\ep

Noticing  that
$\displaystyle{\sum_{i=0}^{\infty}{m\choose i}z^{m-i}a_{i}}$
is an infinite sum  (although it is a finite sum after
applied to each vector $u$ of $M$), we may consider
a certain  completion of $g(V)$. By considering the tensor product
vertex algebra ${\C}((t))\otimes V$ we obtain a Lie algebra
$g_{0}({\C}((t))\otimes V)$ (from Lemma 2.1). It is clear that this
Lie algebra is the completion of $g(V)$ with respect to a certain
topology for $g(V)$. We denote this Lie algebra by $\bar{g}(V)$.

For any
$f(t)=\sum_{m\ge k}c_{m}t^{m}\in {\C}((t))$, since the
following sum:
\begin{eqnarray}
\sum_{m\ge k}c_{m}\left(\sum_{i=0}^{\infty}{m\choose i}z^{m-i}t^{i}\right)
=\sum_{i=0}^{\infty}\left(\sum_{m\ge k}{m\choose i}c_{m}z^{m-i}\right)t^{i}
\end{eqnarray}
may not be a well-defined element of ${\C}((t))$, we cannot extend
an evaluation $g(V)$-module $M_{z}$ to a $\bar{g}(V)$-module. 

Define a linear map $\Delta_{z}$ as follows:
\begin{eqnarray}
\Delta_{z}:& &{\C}[t,t^{-1}]\otimes V\rightarrow ({\C}((t))\otimes
V)\otimes ({\C}((t))\otimes V);\nonumber\\
& &f(t)\otimes a\mapsto 1\otimes (f(t)\otimes a)+
(f(z+t)\otimes a)\otimes 1.
\end{eqnarray}

\begin{prop}\label{p2.12}
$\Delta_{z}$ induces an associative algebra
homomorphism from $U(g(V))$ to $U(\bar{g}(V))\otimes 
U(\bar{g}(V))$.
\end{prop}

{\bf Proof.} Define a linear map $\Delta_{z}^{1}$ from
${\C}[t,t^{-1}]\otimes V$ to ${\C}((t))\otimes V$ as follows:
\begin{eqnarray}
\Delta_{z}^{1}(f(t)\otimes a)=f(z+t)\otimes a\;\;\;\mbox{for }f(t)\in
{\C}[t,t^{-1}],a\in V.
\end{eqnarray}
Then $\Delta_{z}=\Delta_{z}^{1}\otimes 1+1\otimes {\rm id}$. Therefore
it suffices to prove that $\Delta_{z}^{1}$ induces a Lie algebra
homomorphism from $g(V)$ to $\bar{g}(V)$. Let $\psi_{z}$ be the
algebra homomorphism from ${\C}[t,t^{-1}]$ to ${\C}((t))$
defined by: $\psi_{z}(f(t))=f(z+t)$ for $f(t)\in {\C}((t))$. From
Examples 2.4 $\psi_{z}$ is a vertex algebra homomorphism from ${\bf
C}[t,t^{-1}]$ to ${\C}((t))$, so that $\psi_{z}\otimes {\rm id}$ is
a vertex algebra homomorphism from ${\C}[t,t^{-1}]\otimes V$ to ${\bf
C}((t))\otimes V$. By definition $\Delta_{z}^{1}=\psi_{z}\otimes {\rm
id}$. Therefore $\Delta_{z}^{1}$ induces a Lie algebra homomorphism from 
$g(V)$ to $\bar{g}(V)$. $\;\;\;\;\Box$

\begin{rem}\label{r2.13}
 The Hopf-like algebra $(U(g(V)),U(\bar{g}(V)),
\Delta_{z})$ is implicitly used in many references such as
[HL0-4], [KL0-2] and [MS]. 
\end{rem}

\section{A definition of tensor product and a construction}
In this section we shall first formulate a definition of a tensor product in
terms of a certain universal property as an analogue of the notion of
the classical tensor product. Then we give a construction of a tensor product
for an ordered pair of modules for a rational vertex operator algebra.

Throughout this section, $V$ will be a fixed vertex operator algebra.
Let $\cal{C}$ be the category of weak $V$-modules where a morphism $f$ from 
$W$ to $M$ is a linear map such that $f(Y(a,x)u)=Y(a,x)f(u)$ for any 
$a\in V, u\in W.$  Let ${\cal{C}}_{0}$ be the subcategory of $\cal{C}$ where
objects of ${\cal{C}}_{0}$ are weak $V$-modules satisfying all the axioms of
a module except that in (M0), infinite-dimensional homogeneous subspaces 
are allowed. 

\bd{d3.1} 
Let $\cal{D}$ be either the category $\cal{C}$ or ${\cal{C}}_{0}$ and 
let $W_{1}$ and $W_{2}$ be objects of $\cal{D}$. A
{\em tensor product} for the ordered pair $(W_{1},W_{2})$ is a pair
$(M,F(\cdot,x))$ consisting of an object $M$ of $\cal{D}$ and an intertwining
operator $F(\cdot,x)$ of type
$\left(\!\begin{array}{c}M\\W_{1} W_{2}\end{array}\!\right)$
satisfying the following universal 
property: For any object $W$ of $\cal{D}$ and any intertwining
operator $I(\cdot,x)$ of type
$\left(\!\begin{array}{c}W\\W_{1} W_{2}\end{array}\!\right)$, there
exists a unique $V$-homomorphism $\psi$ from $M$ to $W$ such that
$I(\cdot,x)=\psi\circ F(\cdot,x)$. (Here $\psi$ extends canonically to a
linear map from $M\{x\}$ to $W\{x\}$.)
\ed

\begin{rem}\label{r3.2} 
Just as in the classical algebra theory, it follows
{}from the universal property that if there exists a tensor product
$(M,F(\cdot,x))$ in the category $\cal{C}$ or ${\cal{C}}_{0}$
for the ordered pair $(W_{1},W_{2})$, then it is unique up to 
$V$-module isomorphism, i.e., if $(W,G(\cdot,x))$ is another tensor
product, then there is a $V$-module isomorphism $\psi$ from $M$ to $W$
such that $G=\psi \circ F$. Conversely, let $(M,F(\cdot,x))$ be a tensor
product for the ordered pair $(W_{1},W_{2})$ and let $\sigma$ be an
automorphism of the $V$-module $M$. Then $(M, \sigma\circ F(\cdot,x))$
is a tensor product for $(W_{1},W_{2})$.
\end{rem}

\begin{lem}\label{l3.3}
Let $(W,F(\cdot,x))$ is a tensor product in the category 
$\cal{C}$ or ${\cal{C}}_{0}$ for the
ordered pair $(W_{1},W_{2})$. Then $F(\cdot,x)$ is surjective in the
sense that all the coefficients of $F(u_{1},x)u_{2}$ for $u_{i}\in
W_{i}$ linearly span $W$.
\end{lem}

{\bf Proof.} Let $\overline{W}$ be the linear span of all the coefficients
of $F(u_{1},x)u_{2}$ for $u_{i}\in W_{i}$. Then it follows from the commutator formula (\ref{ec}) that
$\overline{W}$ is a
submodule of $W$ and $F(\cdot,x)$ is an intertwining operator
of type
$\left(\!\begin{array}{c}\overline{W}\\W_{1} W_{2}\end{array}\!\right)$. 
It follows from the
universal property of $(W, F(\cdot,x))$ that there is a unique $V$-module
homomorphism $\psi$ from $W$ to $\overline{W}$ such that 
\begin{eqnarray}
F(u_{1},x)u_{2}=\psi(F(u_{1},x)u_{2})\;\;\;\mbox{ for }u_{1}\in W_{1},
u_{2}\in W_{2}.
\end{eqnarray}
Since $\bar{W}$ is a submodule of $W$, $\psi$ may be viewed as a 
$V$-homomorphism from $W$ to $W$.
It follows from the
universal property of $(W, F(\cdot,x))$ that $\psi =1$. Thus $W=\bar{W}$.
Then the proof is complete.$\;\;\;\;\Box$

\begin{coro}\label{c3.4} 
If $(M,F(\cdot,x))$ is a tensor product in the category $\cal{C}$ 
or ${\cal{C}}_{0}$ for the
ordered pair $(W_{1},W_{2})$, then for any weak $V$-module
$W_{3}$ in the same category, ${\rm Hom}_{V}(M,W_{3})$ is naturally
isomorphic to the space
of intertwining operators of type
$\left(\!\begin{array}{c}W_{3}\\W_{1} W_{2}\end{array}\!\right)$.
\end{coro}

{\bf Proof.} Let $\phi$ be any $V$-homomorphism from $M$ to $W_{3}$.
Then $\phi F(\cdot,x)$ is an intertwining operator of type
$\left(\!\begin{array}{c}W_{3}\\W_{1} W_{2}\end{array}\!\right)$. Thus we
obtain a linear map $\pi$ from ${\rm Hom}_{V}(M,W_{3})$ to
$I\left(\!\begin{array}{c}W_{3}\\W_{1} W_{2}\end{array}\!\right)$ defined
by $\pi (\phi)=\phi F(\cdot,x)$. Since $F(\cdot,x)$ is surjective
(Lemma 3.3), $\pi$ is injective. On the other hand, the universal
property of $(W, F(\cdot,x))$ implies that $\pi$ is
surjective.$\;\;\;\;\Box$ 
 
\begin{rem}\label{r3.5}
If $(M,F(\cdot,x))$ is a tensor product in the category $\cal{C}$ 
or ${\cal{C}}_{0}$ for the ordered pair $(W_{1},W_{2})$, 
then one can show that $(M,F^{t}(\cdot,x))$ is a
tensor product in the same category for the ordered pair $(W_{2},W_{1})$. 
This gives a
sort of commutativity of tensor product. It is important to notice that 
it should not be confused with the symmetric property of a tensor category.
As a matter of fact, the tensor category of $V$-modules is not a symmetric
tensor category [HL0-4].
If $(M,Y_{M}(\cdot,x))$ is a $V$-module, one can show that 
$(M,Y_{M}(\cdot,x))$ is a tensor product
for $(V,M)$. This shows that the adjoint module $V$ satisfies a certain unital 
property. 
\end{rem}

Next toward a construction of a tensor product we 
shall construct an ${\N}$-gradable weak $V$-module $T(W_{1},W_{2})$
 for an ordered pair $(W_{1},W_{2}) $ of ${\N}$-gradable weak
$V$-modules. First form the vector space 
\begin{eqnarray}
F_{0}(W_{1},W_{2})={\C}[t,t^{-1}]\otimes W_{1}\otimes W_{2}
\end{eqnarray}
and set
\begin{eqnarray}
Y_{t}(u,x)=\sum_{n\in {\Z}}(t^{n}\otimes u)x^{-n-1}\;\;\;\mbox{for
}u\in W_{1}. 
\end{eqnarray}
Then ${\C}[t,t^{-1}]\otimes W_{1}$ is linearly spanned by the
coefficients of all $Y_{t}(u,x)$ for $u\in W_{1}$.

Fix a gradation $W_{i}=\oplus_{n\in {\N}}W_{i}(n)$ for each $W_{i}$ $(i=1,2)$.
Later we will show that if $V$ is rational, there is a canonical gradation 
for any ${\N}$-gradable weak $V$-module.

We define a ${\Z}$-grading for $F_{0}(W_{1},W_{2})$ as follows:
For $k \in {\Z}; m, n \in {\N},u \in W_{1}(m),v \in W_{2}(n)$, we define
\begin{eqnarray}
\deg\:(t^{k}\otimes u\otimes v)=m+n-k-1.
\end{eqnarray}

Define an action of $\hat{V}={\C}[t,t^{-1}]\otimes V$ on
$F_{0}(W_{1},W_{2})$ as follows: For $a \in V,u \in W_{1},v\in W_{2}$, we 
define
\begin{eqnarray}\label{e3.5}
& &Y_{t}(a,x_{1})(Y_{t}(u,x_{2})\otimes v)\nonumber\\
&=&Y_{t}(u,x_{2})\otimes
Y(a,x_{1})v+{\rm Res}_{x_{0}}x_{2}^{-1}\delta\left(\frac{x_{1}-x_{0}}
{x_{2}}\right)Y_{t}(Y(a,x_{0})u,x_{2})\otimes v.
\end{eqnarray}

\begin{prop}\label{p3.6} Under the above defined action of $\hat{V}$,
$F_{0}(W_{1},W_{2})$ becomes a ${\Z}$-graded $g(V)$-module of
level one, i.e.,
\begin{eqnarray}
& &Y_{t}({\bf 1},x)=1,\;\;Y_{t}(L(-1)a,x)={d \over
dx}Y_{t}(a,x)\;\;\;\mbox{{\it for} }a \in V; \\
& &\deg\:a(n)={\rm wt}\:a-n-1\;\;\;\mbox{{\it for each
homogeneous }}a \in V, n \in {\Z};\\ 
& &[Y_{t}(a,x_{1}),Y_{t}(b,x_{2})]={\rm Res}_{x_{0}}x_{2}^{-1}\delta\left(
\frac{x_{1}-x_{0}}{x_{2}}\right)Y_{t}(Y(a,x_{0})b,x_{2})\;\;\;\mbox{{\it
for }}a,b\in V.
\end{eqnarray}
\end{prop}

{\bf Proof.} Writing (\ref{e3.5}) into components we get
\begin{eqnarray}
& &(t^{m}\otimes a)(t^{n}\otimes u\otimes v)\nonumber\\
&=&t^{n}\otimes u\otimes a_{m}u\nonumber\\
& &+{\rm Res}_{x_{0}}{\rm Res}_{x_{1}}{\rm Res}_{x_{2}}
x_{1}^{m}x_{2}^{n}x_{1}^{-1}x_{1}^{-1}\delta\left(\frac{x_{2}+x_{0}}{x_{1}}
\right)Y_{t}(Y(a,x_{0})u,x_{2})\otimes v  \nonumber\\ 
&=&t^{n}\otimes u\otimes a_{m}u
+{\rm Res}_{x_{2}}\sum_{i=0}^{\infty}
\left(\begin{array}{c}m\\i\end{array}\right)
x_{2}^{m+n-i}Y_{t}(a_{i}u,x_{2})\otimes v\nonumber\\
&=&t^{n}\otimes u\otimes a_{m}u
+\sum_{i=0}^{\infty}\left(\begin{array}{c}m\\i\end{array}\right)
(t^{m+n-i}\otimes a_{i}u\otimes v).
\end{eqnarray}
It follows from Lemma \ref{l2.11} that (\ref{e3.5}) defines a $g(V)$-module
structure on ${\C}[t,t^{-1}]\otimes W_{1}\otimes W_{2}$, which is a
tensor product module of level-zero $g(V)$-module ${\bf
C}[t,t^{-1}]\otimes W_{1}$ with the level-one $g(V)$-module $W_{2}$.
$\;\;\;\;\Box$

Define $J_{0}$ to be the $g(V)$-submodule of $F_{0}(W_{1},W_{2})$ generated 
by the following subspaces:
\begin{eqnarray}
t^{m+n+i} \otimes W_{1}(m) \otimes W_{2}(n)\;\;\;\mbox{for
}m,n,i \in {\N}. 
\end{eqnarray}
Set
\begin{eqnarray}
F_{1}(W_{1},W_{2})=F_{0}(W_{1},W_{2})/J_{0}. 
\end{eqnarray}

\begin{rem}\label{r3.7}
The space $F_{1}(W_{1},W_{2})$ is an ${\N}$-gradable  $g(V)$-module of 
level one, so that the axioms (M1), (M2) and the
commutator formula (\ref{ec}) automatically hold. Furthermore, for any $a \in
V,w \in F_{1}(W_{1},W_{2})$, $Y_{t}(a,x)w$ involves only
finitely many negative powers of $x$.
\end{rem}

\begin{rem}\label{r3.8}
 Notice that the action (\ref{e3.5}) of $g(V)$ on 
$F_{0}(W_{1},W_{2})$ only 
reflects the commutator formula (\ref{ec}), which is weaker than the Jacobi
identity, unlike the situation in the classical Lie
algebra theory. In the next step,  we consider the whole Jacobi identity
relation for an intertwining operator. This step in our approach
might be related to the ``compatibility condition'' in Huang and Lepowsky's
approach [HL0-4]. 
\end{rem}

Let $\pi_{1}$ be the quotient map from $F_{0}(W_{1},W_{2})$ onto
$F_{1}(W_{1},W_{2})$ and
let $J_{1}$ be the subspace of $F_{1}(W_{1},W_{2})$,
linearly spanned by all coefficients of monomials
$x_{0}^{m}x_{1}^{n}x_{2}^{k}$ in the following expressions:
\begin{eqnarray}\label{e3.12}
& &x_{0}^{-1}\delta\left(\frac{x_{1}-x_{2}}{x_{0}}\right)Y_{t}(a,x_{1})\pi
_{1} (Y_{t}(u,x_{2})\otimes v)\nonumber\\
& &-x_{0}^{-1}\delta\left(\frac{x_{2}-x_{1}}{-x_{0}}\right)\pi _{1}
(Y_{t}(u,x_{2})\otimes Y(a,x_{1})v)\nonumber\\
&-&x_{2}^{-1}\delta\left(\frac{x_{1}-x_{0}}{x_{2}}\right)\pi _{1}
(Y_{t}(Y(a,x_{0})u,x_{2}))\otimes v  \end{eqnarray}
for $a \in V, u \in W_{1}, v \in W_{2}$.

\begin{prop}\label{p3.9}
The  subspace $J_{1}$ is a graded $g(V)$-submodule of $F_{1}(W_{1},W_{2})$.
\end{prop}

{\bf Proof}. For $a,b\in V,u\in W_{1},v\in W_{2}$, we have
\begin{eqnarray}
& &x_{0}^{-1}\delta\left(\frac{x_{1}-x_{2}}{x_{0}}\right)Y_{t}(b,x_{3})
Y_{t}(a,x_{1})\pi _{1}(Y_{t}(u,x_{2})\otimes v)\nonumber \\
&=&x_{0}^{-1}\delta\left(\frac{x_{1}-x_{2}}{x_{0}}\right) Y_{t}(a,x_{1})
 Y_{t}(b,x_{3})\pi _{1}(Y_{t}(u,x_{2})\otimes v)\nonumber \\
& &+{\rm Res}_{x_{4}}x_{0}^{-1}\delta\left(\frac{x_{1}-x_{2}}{x_{0}}\right)
x_{1}^{-1}\delta\left(\frac{x_{3}-x_{4}}{x_{1}}\right)Y_{t}(Y(b,x_{4})a,x_{1})
\pi _{1}(Y_{t}(u,x_{2})\otimes v)\nonumber \\
&=&x_{0}^{-1}\delta\left(\frac{x_{1}-x_{2}}{x_{0}}\right)Y_{t}(a,x_{1})
\pi _{1}(Y_{t}(u,x_{2})\otimes Y(b,x_{3})v)\nonumber \\
& &+{\rm Res}_{x_{4}}x_{0}^{-1}\delta\left(\frac{x_{1}-x_{2}}{x_{0}}\right)
x_{2}^{-1}\delta\left(\frac{x_{3}-x_{4}}{x_{2}}\right)Y_{t}(a,x_{1})
\pi_{1}(Y_{t}(Y(b,x_{4})u,x_{2})\otimes v)\nonumber \\
& &+{\rm Res}_{x_{4}}x_{0}^{-1}\delta\left(\frac{x_{1}-x_{2}}{x_{0}}\right)
x_{1}^{-1}\delta\left(\frac{x_{3}-x_{4}}{x_{1}}\right)Y_{t}(Y(b,x_{4})a,x_{1})
\pi_{1}(Y_{t}(u,x_{2})\otimes v); \nonumber\\
&&\mbox{}
\end{eqnarray}
\begin{eqnarray}
& &Y_{t}(b,x_{3})\pi_{1}(Y_{t}(u,x_{2})\otimes Y(a,x_{1})v)\nonumber \\
&=&\pi_{1}(Y_{t}(u,x_{2})\otimes Y(b,x_{3})Y(a,x_{1})v)\nonumber\\
& &+{\rm Res}_{x_{4}}x_{2}^{-1}\delta\left(\frac{x_{3}-x_{4}}{x_{2}}
\right)\pi_{1}
(Y_{t}(Y(b,x_{4})u,x_{2})\otimes Y(a,x_{1})v)\nonumber\\
&=&\pi_{1}(Y_{t}(u,x_{2})\otimes Y(a,x_{1})Y(b,x_{3})v)\nonumber\\
& &+{\rm Res}_{x_{4}}x_{1}^{-1}\delta\left(\frac{x_{3}-x_{4}}{x_{1}}\right)
\pi_{1}
(Y_{t}(u,x_{2})\otimes Y(Y(b,x_{4})a,x_{1})v)\nonumber \\
& &+{\rm Res}_{x_{4}}x_{2}^{-1}\delta\left(\frac{x_{3}-x_{4}}{x_{2}}
\right)\pi_{1}
(Y_{t}(Y(b,x_{4})u,x_{2})\otimes Y(a,x_{1})v); 
\end{eqnarray}
and
\begin{eqnarray}
& &x_{2}^{-1}\delta\left(\frac{x_{1}-x_{0}}{x_{2}}\right)Y_{t}(b,x_{3})
\pi_{1}
(Y_{t}(Y(a,x_{0})u,x_{2})\otimes v)\nonumber\\
&=&x_{2}^{-1}\delta\left(\frac{x_{1}-x_{0}}{x_{2}}\right)
\pi_{1}(Y_{t}(Y(a,x_{0})u,x_{2})
\otimes Y(b,x_{4})v)\nonumber\\
& &+{\rm Res}_{x_{4}}x_{2}^{-1}\delta\left(\frac{x_{1}-x_{0}}{x_{2}}\right)
x_{2}^{-1}\delta\left(\frac{x_{3}-x_{4}}{x_{2}}\right)\pi_{1}(Y_{t}
(Y(b,x_{4})Y(a,x_{0})u,x_{2})\otimes v)\nonumber \\
&=&x_{2}^{-1}\delta\left(\frac{x_{1}-x_{0}}{x_{2}}\right)Y_{t}(b,x_{3})
\pi_{1}
(Y_{t}(Y(a,x_{0})u,x_{2})\otimes v)\nonumber\\
& &+{\rm Res}_{x_{4}}x_{2}^{-1}\delta\left(\frac{x_{1}-x_{0}}{x_{2}}\right)
x_{2}^{-1}\delta\left(\frac{x_{3}-x_{4}}{x_{2}}\right)\pi_{1}(Y_{t}
(Y(a,x_{0})Y(b,x_{4})u,x_{2})\otimes v)\nonumber \\
& &+{\rm Res}_{x_{4}}{\rm Res}_{x_{5}}x_{2}^{-1}\delta\left(
\frac{x_{1}-x_{0}}{x_{2}}\right)
x_{2}^{-1}\delta\left(\frac{x_{3}-x_{4}}{x_{2}}\right) x_{0}^{-1}
\delta\left(\frac{x_{4}-x_{5}}
{x_{0}}\right)\cdot \nonumber \\
& &\cdot \pi_{1}(Y_{t}(Y(Y(b,x_{5})a,x_{0})u,x_{2})\otimes v)\nonumber \\
&=&x_{2}^{-1}\delta\left(\frac{x_{1}-x_{0}}{x_{2}}\right)Y_{t}(b,x_{3})
\pi_{1}
(Y_{t}(Y(a,x_{0})u,x_{2})\otimes v)\nonumber\\
& &+{\rm Res}_{x_{4}}x_{2}^{-1}\delta\left(\frac{x_{1}-x_{0}}{x_{2}}\right)
x_{2}^{-1}\delta\left(\frac{x_{3}-x_{4}}{x_{2}}\right)
x_{2}^{-1}\delta\left(\frac{x_{3}-x_{4}}{x_{2}}\right)\nonumber \\
& &\cdot\pi_{1}(Y_{t}
(Y(a,x_{0})Y(b,x_{4})u,x_{2})\otimes v)\nonumber \\
& &+{\rm Res}_{x_{5}}x_{2}^{-1}\delta\left(\frac{x_{1}-x_{0}}{x_{2}}\right)
x_{2}^{-1}\delta\left(\frac{x_{3}-x_{0}-x_{5}}{x_{2}}\right)\cdot
\nonumber \\
& &\cdot \pi_{1}(Y_{t}(Y(Y(b,x_{5})a,x_{0})u,x_{2})\otimes v)\nonumber \\
&=&x_{2}^{-1}\delta\left(\frac{x_{1}-x_{0}}{x_{2}}\right)Y_{t}(b,x_{3})
\pi_{1}
(Y_{t}(Y(a,x_{0})u,x_{2})\otimes v)\nonumber\\
& &+{\rm Res}_{x_{4}}x_{2}^{-1}\delta\left(\frac{x_{1}-x_{0}}{x_{2}}\right)
x_{2}^{-1}\delta\left(\frac{x_{3}-x_{4}}{x_{2}}\right)
x_{2}^{-1}\delta\left(\frac{x_{3}-x_{4}}{x_{2}}\right)\nonumber \\
& &\cdot \pi_{1}(Y_{t}
(Y(a,x_{0})Y(b,x_{4})u,x_{2})\otimes v)\nonumber \\
& &+{\rm Res}_{x_{4}}x_{2}^{-1}\delta\left(\frac{x_{1}-x_{0}}{x_{2}}
\right)x_{1}^{-1}
\delta\left(\frac{x_{3}-x_{4}}{x_{1}}\right)\pi_{1}
(Y_{t}(Y(Y(b,x_{5})a,x_{0})u,x_{2}) \otimes v).\nonumber\\
& &\mbox{}
\end{eqnarray}
Then it is clear that $J_{1}$ is stable under the action of 
$Y_{t}(b,x)$ for any $b\in V$.  $\;\;\;\;\Box$

\bt{t3.10}
The quotient space $F_{2}(W_{1},W_{2}):=F_{1}(W_{1},W_{2})/J_{1}$ is an
${\N}$-gradable weak $V$-module.
\et

{\bf Proof}. We only need to prove the Jacobi identity. Let $\pi _{2}$ be the natural quotient map from $F_{0}(W_{1},W_{2})$ onto $F_{2}(W_{1},W_{2})$.
For $a,b\in
V,u\in W_{1},v\in W_{2}$, we have
\begin{eqnarray}
& &x_{0}^{-1}\delta\left(\frac{x_{1}-x_{2}}{x_{0}}\right)Y_{t}(a,x_{1})
Y_{t}(b,x_{2})\pi_{2}(Y_{t}(u,x_{3})\otimes v)   \nonumber \\
&=&x_{0}^{-1}\delta\left(\frac{x_{1}-x_{2}}{x_{0}}\right)Y_{t}(a,x_{1})
\pi_{2}\left( Y_{t}(u,x_{3})\otimes Y(b,x_{2})v\right) \nonumber \\
& &+{\rm Res}_{x_{4}}x_{0}^{-1}\delta\left(\frac{x_{1}-x_{2}}{x_{0}}\right)
x_{3}^{-1}\delta\left(\frac{x_{2}-x_{4}}{x_{3}}\right)Y_{t}(a,x_{1})\pi_{2}
\left(
Y_{t}(Y(b,x_{4})u,x_{3})\otimes v\right)\nonumber \\
&=&x_{0}^{-1}\delta\left(\frac{x_{1}-x_{2}}{x_{0}}\right)\pi_{2}(Y_{t}(u,x_{3})\otimes
Y(a,x_{1})Y(b,x_{2})v)\label{e3.24} \\
& &+ {\rm Res}_{x_{4}}x_{0}^{-1}\delta\left(\frac{x_{1}-x_{2}}{x_{0}}\right)
x_{3}^{-1}\delta\left(\frac{x_{1}-x_{4}}{x_{3}}\right)\pi_{2}(Y_{t}(Y(a,x_{4})u,
x_{3})\otimes Y(b,x_{2})v) \label{e3.25}\nonumber\\
& &\mbox{}\\
& &+ {\rm Res}_{x_{4}}x_{0}^{-1}\delta\left(\frac{x_{1}-x_{2}}{x_{0}}\right)
x_{3}^{-1}\delta\left(\frac{x_{2}-x_{4}}{x_{3}}\right)
Y_{t}(a,x_{1})\pi_{2}(Y_{t}(
Y(b,x_{4})u,x_{3})\otimes v );\label{e3.26}\nonumber\\
& &\mbox{}
\end{eqnarray}
\begin{eqnarray}
& &x_{0}^{-1}\delta\left(\frac{-x_{2}+x_{1}}{x_{0}}\right)Y_{t}(b,x_{2})
Y_{t}(a,x_{1})\pi_{2}(Y_{t}(u,x_{3})\otimes v)   \nonumber \\
&=&x_{0}^{-1}\delta\left(\frac{-x_{2}+x_{1}}{x_{0}}\right)\pi_{2}(Y_{t}(u,x_{3})\otimes
Y(b,x_{2})Y(a,x_{1})v)\label{e3.27}\\
& &+{\rm Res}_{x_{4}}x_{0}^{-1}\delta\left(\frac{-x_{2}+x_{1}}{x_{0}}\right)
x_{3}^{-1}\delta\left(\frac{x_{2}-x_{4}}{x_{3}}\right)\pi_{2}(Y_{t}(Y(b,x_{4})u,
x_{3})\otimes Y(a,x_{1})v)\label{e3.28}\nonumber\\
& &\mbox{}\\
& &+{\rm Res}_{x_{4}}x_{0}^{-1}\delta\left(\frac{-x_{2}+x_{1}}{x_{0}}\right)
x_{3}^{-1}\delta\left(\frac{x_{1}-x_{4}}{x_{3}}\right)Y_{t}(b,x_{2})\pi_{2}(Y_{t}(
Y(a,x_{4})u,x_{3})\otimes v);\label{e3.29}\nonumber\\
& &\mbox{}
\end{eqnarray}
and
\begin{eqnarray}
& &x_{2}^{-1}\delta\left(\frac{x_{1}-x_{0}}{x_{2}}\right)Y_{t}(Y(a,x_{0})
b,x_{2})\pi_{2}(Y_{t}(u,x_{3})\otimes v)   \nonumber \\
&=&x_{2}^{-1}\delta\left(\frac{x_{1}-x_{0}}{x_{2}}\right)\pi_{2}(Y_{t}(u,x_{3})\otimes
Y(Y(a,x_{0})b,x_{2})v)\label{e3.30}\\
& &+{\rm Res}_{x_{4}}x_{2}^{-1}\delta\left(\frac{x_{1}-x_{0}}{x_{2}}\right)
x_{3}^{-1}\delta\left(\frac{x_{2}-x_{4}}{x_{3}}\right)\pi_{2}(Y_{t}(Y(Y(a,x_{0})b,
x_{4})u,x_{3})\otimes v).\label{e3.31}\nonumber\\
& &\mbox{}
\end{eqnarray}
It follows from the Jacobi identity of $W_{2}$ that
$(\ref{e3.24})-(\ref{e3.27})=(\ref{e3.30})$. Since 
\begin{eqnarray}
& &x_{0}^{-1}\delta\left(\frac{x_{1}-x_{2}}{x_{0}}\right)x_{3}^{-1}\delta
\left(\frac{x_{2}-x_{4}}{x_{3}}\right)\nonumber \\
&=&x_{0}^{-1}\delta\left(\frac{x_{1}-x_{3}-x_{4}}{x_{0}}\right)x_{3}^{-1}
\delta\left(\frac{x_{2}-x_{4}}{x_{3}}\right)\nonumber \\
&=&(x_{1}-x_{3})^{-1}\delta\left(\frac{x_{0}+x_{4}}{x_{1}-x_{3}}\right)
x_{3}^{-1}\delta\left(\frac{x_{2}-x_{4}}{x_{3}}\right)\nonumber\\
&=&(x_{0}+x_{4})^{-1}\delta\left(\frac{x_{1}-x_{3}}{x_{0}+x_{4}}\right)
x_{3}^{-1}\delta\left(\frac{x_{2}-x_{4}}{x_{3}}\right);
\end{eqnarray}
\begin{eqnarray}
& &x_{0}^{-1}\delta\left(\frac{-x_{2}+x_{1}}{x_{0}}\right)x_{3}^{-1}\delta
\left(\frac{x_{2}-x_{4}}{x_{3}}\right)\nonumber \\
&=&x_{0}^{-1}\delta\left(\frac{-x_{3}-x_{4}+x_{1}}{x_{0}}\right)
x_{3}^{-1}\delta
\left(\frac{x_{2}-x_{4}}{x_{3}}\right)\nonumber \\
&=&(x_{0}+x_{4})^{-1}\delta\left(\frac{-x_{3}+x_{1}}{x_{0}+x_{4}}\right)
x_{3}^{-1}\delta\left(\frac{x_{2}-x_{4}}{x_{3}}\right),
\end{eqnarray}
by the $J_{1}$-defining relation (\ref{e3.12}), we have
\begin{eqnarray}
& &(\ref{e3.26})-(\ref{e3.28})\nonumber \\
&=&{\rm Res}_{x_{4}}x_{3}^{-1}\delta
\left(\frac{x_{2}-x_{4}}{x_{3}}\right)x_{3}^{-1}\delta\left(\frac{x_{1}
-x_{0}-x_{4}}{x_{3}}\right)\nonumber \\
& &\cdot \pi_{2}(Y_{t}(
Y(a,x_{0}+x_{4})Y(b,x_{4})u,x_{3})\otimes v)\nonumber\\
&=&{\rm Res}_{x_{4}}{\rm Res}_{x_{5}}x_{5}^{-1}\delta\left(\frac{x_{0}+x_{4}}
{x_{5}}\right)x_{3}^{-1}\delta
\left(\frac{x_{2}-x_{4}}{x_{3}}\right)x_{3}^{-1}\delta\left(\frac
{x_{1}-x_{5}}{x_{3}}\right)\nonumber\\
& &\cdot \pi_{2}(Y_{t}(Y(a,x_{5})Y(b,x_{4})u,x_{3})\otimes v)\nonumber\\
&=&{\rm Res}_{x_{4}}{\rm Res}_{x_{5}}x_{0}^{-1}\delta\left(\frac{x_{5}-x_{4}}
{x_{0}}\right)x_{3}^{-1}\delta
\left(\frac{x_{2}-x_{4}}{x_{3}}\right)x_{3}^{-1}\delta\left(\frac
{x_{1}-x_{5}}{x_{3}}\right)\nonumber\\
& &\cdot \pi_{2}(Y_{t}(Y(a,x_{5})Y(b,x_{4})u,x_{3})\otimes v).
\end{eqnarray}
Similarly, we obtain
\begin{eqnarray}
& &(\ref{e3.25})-(\ref{e3.29})\nonumber \\
&=&-{\rm Res}_{x_{4}}x_{3}^{-1}\delta
\left(\frac{x_{1}-x_{4}}{x_{3}}\right)x_{3}^{-1}\delta\left(\frac
{x_{1}+x_{0}-x_{4}}{x_{3}}\right)\nonumber\\
& &\cdot \pi_{2}(Y_{t}(Y(a,-x_{0}+x_{4})Y(b,x_{4})u,x_{3})\otimes v)\nonumber\\
&=&-{\rm Res}_{x_{4}}{\rm Res}_{x_{5}}x_{5}^{-1}\delta\left(\frac{-x_{0}+x_{4}}
{x_{5}}\right)x_{3}^{-1}\delta
\left(\frac{x_{1}-x_{4}}{x_{3}}\right)x_{3}^{-1}\delta\left(\frac{x_{2}
-x_{5}}{x_{3}}\right)\nonumber \\
& &\cdot \pi_{2}(Y_{t}(Y(b,x_{4})Y(a,x_{5})u,x_{3})\otimes v)\nonumber\\
&=&-{\rm Res}_{x_{4}}{\rm Res}_{x_{5}}x_{0}^{-1}\delta\left(\frac{-x_{5}+x_{4}}
{x_{0}}\right)x_{3}^{-1}\delta
\left(\frac{x_{1}-x_{4}}{x_{3}}\right)x_{3}^{-1}\delta\left(\frac{x_{2}
-x_{5}}{x_{3}}\right)\nonumber \\
& &\cdot \pi_{2}(Y_{t}(Y(b,x_{4})Y(a,x_{5})u,x_{3})\otimes v) \nonumber\\
&=&-{\rm Res}_{x_{5}}{\rm Res}_{x_{4}}x_{0}^{-1}\delta\left(\frac{-x_{4}+x_{5}}
{x_{0}}\right)x_{3}^{-1}\delta
\left(\frac{x_{1}-x_{5}}{x_{3}}\right)x_{3}^{-1}\delta\left(\frac{x_{2}
-x_{4}}{x_{3}}\right)\nonumber \\
& &\cdot \pi_{2}(Y_{t}(Y(b,x_{5})Y(a,x_{4})u,x_{3})\otimes v).
\end{eqnarray}
Therefore, we have
\begin{eqnarray}
& &(\ref{e3.25})+(\ref{e3.26})-(\ref{e3.28})-(\ref{e3.29})\nonumber \\
&=&{\rm Res}_{x_{4}}{\rm Res}_{x_{5}}x_{4}^{-1}\delta\left(
\frac{x_{5}-x_{0}}{x_{4}}\right)x_{3}^{-1}\delta
\left(\frac{x_{2}-x_{4}}{x_{3}}\right)
x_{3}^{-1}\delta\left(\frac{x_{1}-x_{5}}{x_{3}}\right)\nonumber \\
& &\cdot \pi_{2}(Y_{t}(Y(Y(a,x_{0})b,x_{4})u,x_{3})\otimes v)\nonumber \\
&=&{\rm Res}_{x_{4}}x_{3}^{-1}\delta
\left(\frac{x_{2}-x_{4}}{x_{3}}\right)x_{2}^{-1}\delta\left(\frac{x_{1}-x_{0}}
{x_{2}}\right)\pi_{2}(Y_{t}(Y(Y(a,x_{0})b,x_{4})u,x_{3})\otimes v)\nonumber\\
&=& (\ref{e3.31}).\end{eqnarray}
Here we have used the following fact:
\begin{eqnarray}
& &{\rm Res}_{x_{5}}x_{4}^{-1}\delta\left(
\frac{x_{5}-x_{0}}{x_{4}}\right)x_{3}^{-1}\delta
\left(\frac{x_{2}-x_{4}}{x_{3}}\right)
x_{3}^{-1}\delta\left(\frac{x_{1}-x_{5}}{x_{3}}\right)\nonumber \\
&=&{\rm Res}_{x_{5}}x_{3}^{-1}\delta
\left(\frac{x_{2}-x_{4}}{x_{3}}\right)x_{3}^{-1}\delta\left(\frac{x_{1}
-x_{4}-x_{0}}{x_{3}}\right)x_{5}^{-1}\delta\left(\frac{x_{4}+x_{0}}{x_{5}}
\right) \nonumber \\
&=&x_{3}^{-1}\delta
\left(\frac{x_{2}-x_{4}}{x_{3}}\right)x_{1}^{-1}\delta\left(\frac{x_{2}-x_{4}
+x_{4}+x_{0}}{x_{1}}\right)\nonumber \\
&=&x_{3}^{-1}\delta
\left(\frac{x_{2}-x_{4}}{x_{3}}\right)x_{2}^{-1}\delta\left(\frac{x_{1}
-x_{0}}{x_{2}}\right). \end{eqnarray}
Then the Jacobi identity is proved. $\;\;\;\;\Box$

Since $F_{2}(W_{1},W_{2})$ is a weak $V$-module, we will freely use
$Y(a,x)$ for $Y_{t}(a,x)$. Recall from Section 2 that 
any weak $V$-module $M$ is a
restricted $g(V)$-module and that for any restricted $g(V)$-module $M$, 
$\bar{M}:=M/J(M)$ 
is a weak $V$-module, where $J(M)$ is the intersection of all ker $f$ with
 $f$ running through 
all $g(V)$-homomorphisms from $M$ to weak $V$-modules. Then 
Theorem \ref{t3.10} says that $J_{1}= J(F_{1}(W_{1},W_{2}))$.

\br{rgi}
If $W_{i}$ for $i=1,2,3$ are just restricted $g(V)$-modules, we can also define
intertwining operator by using the same axioms as in Definition \ref{d2.5}. 
Then following the proof given in [FHL] one can easily see that the 
transpose operator $I^{t}(\cdot,x)$
is well defined and it is an intertwining operator.
\er

\bp{p3.11}
Let $W_{1}$ and $W_{3}$ be weak
$V$-modules, let $M$ be a restricted $g(V)$-module and
let $I(\cdot,x)\!$ be an intertwining 
operator of type 
$\!\left(\!\!\begin{array}{c}W_{3}\\W_{1} M\end{array}\!\!\right)\!.\!$ Then
$I(\cdot,x)J(M)\!=0$, so that we obtain an intertwining operator of type
 $\left(\!\begin{array}{c}W_{3}\\W_{1} \bar{M}\end{array}\!\right)$.
\ep

{\bf Proof.} In the proof of Theorem \ref{t3.10}, replace $W_{2}$,
$F_{2}(W_{1},W_{2})$ and $Y_{t}(\cdot,x)$
 by $M$,  $W_{3}$ and $I(\cdot,x)$, respectively.
Then the $J_{1}$-defining relation (\ref{e3.12}) or Jacobi 
identity for $I(\cdot,x)$ and the Jacobi identity for $W_{3}$:
$$(\ref{e3.25})+(\ref{e3.26})-(\ref{e3.28})-(\ref{e3.29})= (\ref{e3.31})$$
are given.
Following the proof of Theorem \ref{t3.10}, we obtain
$(\ref{e3.24})-(\ref{e3.27})=(\ref{e3.30})$. That is, $I(\cdot,x)J(M)=0$.
Then the proof is complete.
$\;\;\;\;\Box$

Symmetrically, we have

\bp{p3.12} 
Let $W_{2}$ and $W_{3}$ be 
weak $V$-modules, let $M$ be a restricted $g(V)$-module and
let $I(\cdot,x)$ be an
intertwining operator of type 
$\left(\!\begin{array}{c}W_{3}\\M W_{2}\end{array}\!\right)$. Then
$I(J(M),x)=0$, so that we obtain an intertwining operator of type
 $\left(\!\begin{array}{c}W_{3}\\\bar{M} W_{2}\end{array}\!\right)$.
\ep

{\bf Proof.} The proof of this proposition does not directly follow
{}from the proof of Theorem \ref{t3.10}, but it follows from Proposition
\ref{p3.11} and the notion of transpose intertwining operator. 
Since
$I^{t}(\cdot,x)$ is an intertwining operator of type
$\left(\!\begin{array}{c}W_{3}\\W_{2} M\end{array}\!\right)$, by
Proposition \ref{p3.11} we get $I^{t}(\cdot,x)J(M)=0$. Thus $I(J(M),x)=0$.
$\;\;\;\;\Box$

To construct a tensor product out of the weak $V$-module $F_{2}(W_{1},W_{2})$,
we shall study the universal property.
{\em For simplicity, from now on we assume that $W_{1}$ and $W_{2}$ are
weak $V$-modules in the category ${\cal{C}}_{0}$ such that
$W_{i}=\oplus_{n\in {\N}}(W_{i})_{(n+h_{i})}$ for $i=1,2.$}

Let $W$ be a weak $V$-module in the category ${\cal{C}}_{0}$ such that 
$W=\oplus_{n=0}^{\infty}W_{(n+h)}$ for some
 $h$. Let $I(\cdot,x)$ be an
intertwining operator of type
$\left(\!\begin{array}{c}W\\W_{1}W_{2}\end{array}\!\right)$. 
Let $I^{o}(\cdot,x)=x^{h_{1}+h_{2}-h}I(\cdot,x)$ be the normalization.
Then we define
\begin{eqnarray}
\psi_{I}: F_{0}(W_{1},W_{2})\rightarrow W,\;t^{n}\otimes
u \otimes v \mapsto I_{u}(n)v \end{eqnarray}
for $u \in W_{1},v \in W_{2},n \in {\Z}.$
In terms of generating elements, $\psi_{I}$ can be written as:
\begin{eqnarray}
\psi_{I}(Y_{t}(u,x)\otimes v)=
I^{o}(u,x)v\;\;\;\mbox{for }u \in W_{1},v \in W_{2}.
\end{eqnarray}

\begin{lem}\label{l3.13}
The linear map $\psi_{I}$ is a $g(V)$-homomorphism. In other words,
\begin{eqnarray}\label{e3.35}
\psi_{I}(Y_{t}(a,x)w)=Y(a,x)\psi_{I}(w)\;\;\;\mbox{{\it for} } a \in
V, w \in F_{0}(W_{1},W_{2}). 
\end{eqnarray}
\end{lem}

{\bf Proof.} For $a \in V, u \in W_{1},v \in W_{2}$, we have
\begin{eqnarray*}& &\psi_{I}(Y_{t}(a,x_{1})(Y_{t}(u,x_{2})\otimes v))\\
&=&\psi_{I}(Y_{t}(u,x_{2})\otimes Y(a,x_{1})v)
+{\rm Res}_{x_{0}}x_{2}^{-1}\delta\left(\frac{x_{1}-x_{0}}{x_{2}}\right)\psi_{I}(Y_{t}(Y(a,x_{0})u,x_{2})\otimes
v)\\ 
&=&I^{o}(u,x_{2})Y(a,x_{1})v+{\rm Res}_{x_{0}}x_{2}^{-1}\delta\left(\frac{x_{1}-x_{0}}{x_{2}}\right)I^{o}(Y(a,x_{0})u,x_{2})v\\
&=&Y(a,x_{1})I^{o}(u,x_{2})v \\
&=&Y(a,x_{1})\psi_{I}(Y_{t}(u,x_{2})\otimes v).\;\;\;\;\Box
\end{eqnarray*}

\begin{coro}\label{c3.14}
The linear map $\psi_{I}$ induces a
$V$-homomorphism $\bar{\psi}_{I}$ from
$F_{2}(W_{1},W_{2})$ to $W$ such that $\bar{\psi}_{I}$
preserves the 
${\N}$-gradation and $\bar{\psi}_{I}=\pi_{2} \psi_{I}$, where
$\pi_{2}$ is the quotient map from $F_{0}(W_{1},W_{2})$ to
$F_{2}(W_{1},W_{2})$. 
\end{coro}

{\bf Proof}. From Proposition \ref{p2.6} and the Jacobi identity for a 
$V$-module and for an intertwining
operator we get: $J_{0}+J_{1} \subseteq {\rm ker}\;\psi _{I}$. Then we
 have an 
induced linear map $\bar{\psi}_{I}$ from
$F_{2}(W_{1},W_{2})$ to $W$. From (\ref{e3.35})
$\bar{\psi}_{I}$ is a $V$-homomorphism. $\;\;\;\;\Box$

Let $W=\oplus _{n\in {\N}}W_{(n+h)}$ be given as before.
Let $\mbox{Hom}^{0}_{V}(F_{2}(W_{1},W_{2}),W)$ be the space of all
$V$-homomorphisms from $F_{2}(W_{1},W_{2})$ to $W$ which
preserve the given ${\N}$-gradation.
Then we define the following linear map:
\begin{eqnarray}
\bar{\psi}:& &I\!\left(\!\begin{array}{c}W\\W_{1}W_{2}\end{array}\!\right)
 \rightarrow  
\mbox{Hom}_{V}^{0}(F_{2}(W_{1},W_{2}),W)\nonumber\\
& &I(\cdot,x)\mapsto
\bar{\psi}_{I},
\end{eqnarray}

\bp{p3.15}
The map $\bar{\psi} :
I\!\left(\!\begin{array}{c}W\\W_{1}W_{2}\end{array}\!\right)\rightarrow 
\mbox{Hom}^{0}_{V}(F_{2}(W_{1},W_{2}),W); I \mapsto
\bar{\psi}_{I}$ is a linear isomorphism.
\ep

{\bf Proof}. It is clear that $\bar{\psi}$ is injective. Let $f$ be a
$V$-homomorphism from $F_{2}(W_{1},W_{2})$ to $W$ that
preserves the ${\Z}$-grading.
Define a linear map $I(\cdot,x)$ from $W_{1}$ to ${\rm Hom}(W_{2},W)\{x\}$ as 
follows:
\begin{eqnarray}
I(u_{1},x)u_{2}=x^{h-h_{1}-h_{2}}f\pi_{2}(Y_{t}(u_{1},x)\otimes u_{2})
\end{eqnarray}
for any $u_{i}\in W_{i}$. It follows from the defining relations $J_{0}$ and 
$J_{1}$ that  $I(\cdot,x)$ satisfies the axioms (I1) and (I3).
If we prove (I2), then $I(\cdot,x)$ is an intertwining operator satisfying
$\bar{\psi}_{I}=f$. 
For $k \in {\Z},m,n \in {\N}; u \in W_{1}(m), v \in W_{2}$,
we have
 $$\deg\:(t^{k}\otimes u \otimes v)=m+n-k-1.$$
Therefore
\begin{eqnarray}
L(0)f\pi (t^{k}\otimes u\otimes v)=(h+m+n-k-1)f\pi
(t^{k}\otimes u\otimes v). 
\end{eqnarray}
By formula (\ref{e3.5}), we obtain
\begin{eqnarray}
& &L(0)(t^{k}\otimes u \otimes v)\nonumber \\
&=&t^{k}\otimes u \otimes
L(0)v+t^{k+1}\otimes u \otimes v+t^{k}\otimes L(0)u \otimes
v\nonumber \\
&=&t^{k+1}\otimes L(-1)u \otimes v+(h_{1}+h_{2}+m+n)t^{k}\otimes
u\otimes v. 
\end{eqnarray}
Therefore
\begin{eqnarray}\label{e3.36}
f\pi \left(t^{k+1}\otimes L(-1)u\otimes
v+(h_{1}+h_{2}-h+k+1)t^{k}\otimes u\otimes v\right)=0.
\end{eqnarray}
This is exactly the axiom (I2) in terms of components. Then the 
proof is complete.$\;\;\;\;\Box$

For any nonzero ${\N}$-gradable weak $V$-module with a fixed a gradation
$M=\oplus_{n\in {\N}}M(n)$ such that $M(0)\ne 0$, we define 
the {\it radical} of $M$ to be the maximal graded 
submodule $R(M)$ 
such that $R(M)\cap M(0)=0$. 

\bd{d3.16}
 Define $T(W_{1},W_{2})$ to be the
quotient module of $F_{2}(W_{1},W_{2})$ modulo the radical of
$F_{2}(W_{1},W_{2})$ with respect to the given gradation. 
\ed

As a corollary of Proposition \ref{p3.15} we get

\bc{c3.16}The linear isomorphism $\bar{\psi} :
I\!\left(\!\begin{array}{c}W\\W_{1}W_{2}\end{array}\!\right)\rightarrow 
\mbox{Hom}^{0}_{V}(F_{2}(W_{1},W_{2}),W); I \mapsto
\bar{\psi}_{I}$ gives rise to a linear isomorphism from
$I\!\left(\!\begin{array}{c}W\\W_{1}W_{2}\end{array}\!\right) $
to $\mbox{Hom}^{0}_{V}(T(W_{1},W_{2}),W)$, the space of all $V$-homomorphisms 
which preserve the given gradation.
\ec

{}{\em From now on we shall assume that $V$ is rational.}
Then up to equivalence, $V$ has only finitely many irreducible modules. Let 
$\lambda_{1},\cdots, \lambda_{k}$ be all the distinct lowest weights of 
irreducible 
$V$-modules. For any ${\N}$-gradable weak  $V$-module $W$, let $W^{(i)}$ be 
the sum of all irreducible submodules of $W$ with lowest weight $\lambda_{i}$.
Then we obtain a canonical decomposition $W=\oplus_{i=1}^{k}W^{(i)}$. 
Since $W^{(i)}$ is a direct sum of irreducible $V$-modules with lowest 
weight $\lambda_{i}$, any submodule of $W^{(i)}$ is a direct sum of 
irreducible modules with lowest weight $\lambda_{i}$. Thus
$W^{(i)}=\oplus_{n\in {\N}}W^{(i)}_{(n+\lambda_{i})}:
=\oplus_{n\in {\N}}W^{(i)}(n)$. Then 
$\displaystyle{W=\oplus_{n\in {\N}}\oplus_{i=1}^{k}W^{(i)}_{(n+\lambda_{i})}}$.
Then $W$ has a canonical ${\N}$-gradation with 
$\displaystyle{W(n)=\oplus_{i=1}^{k}W^{(i)}_{(n+\lambda_{i})}}$ for 
$n\in {\N}$. It is clear that
 the radical of $W$ is zero with respect to this gradation.
In particular,
$$T(W_{1},W_{2})=T(W_{1},W_{2})^{(1)}\oplus \cdots \oplus 
T(W_{1},W_{2})^{(k)}.$$
Since $F_{2}(W_{1},W_{2})$ is completely reducible, $T(W_{1},
W_{2})$ is isomorphic to the submodule generated by the degree-zero
subspace $F_{2}(W_{1},W_{2})(0)$. 
Let $P_{i}$ be the projection map of $T(W_{1},W_{2})$
onto $T(W_{1},W_{2})^{(i)}$ and let $\pi$ be the natural quotient map from 
$F_{0}(W_{1},W_{2})$ onto 
$T(W_{1},W_{2})$. Then we define
\begin{eqnarray}
F(\cdot,x) :& &W_{1} \rightarrow \left({\rm Hom}_{{\bf
C}}(W_{2},T(W_{1},W_{2})) \right)\{x\};\nonumber\\
& &u_{1}\mapsto F(u_{1},x)\;\;\;\;\mbox{  for }u_{1}\in W_{1}
\end{eqnarray}
where $F(u_{1},x)u_{2}=\sum_{i=1}^{k}x^{\lambda_{i}-h_{1}-h_{2}}P_{i}\pi
(Y_{t}(u_{1},x)\otimes u_{2})$ for $u_{1}\in W_{1}, u_{2}\in W_{2}$.

\bp{p3.17} 
Suppose that $V$ is rational. Then
the defined map $F(\cdot,x)$ is an intertwining operator of type
$\left(\!\begin{array}{c}T(W_{1},W_{2})\\W_{1} W_{2}
\end{array}\!\right)$.
\ep

{\bf Proof.} Let $F_{i}(u_{1},x)u_{2}=x^{\lambda_{i}-h_{1}-h_{2}}P_{i}\pi
(Y_{t}(u_{1},x)\otimes u_{2})$. Then it follows from Proposition \ref{p3.15}
that each $F_{i}(\cdot,x)$ is an intertwining operator of type
$\left(\begin{array}{c}T(W_{1},W_{2})^{(i)}\\W_{1} W_{2}\end{array}\right)$. 
Then it follows immediately. $\;\;\;\;\Box$

\bt{t3.18}
If $V$ is rational and
$W_{i}$ (i=1,2,3) are irreducible weak $V$-modules in the category 
${\cal{C}}_{0}$, then
the pair $(T(W_{1},W_{2}), F(\cdot,x))$ 
is a tensor product in the category ${\cal{C}}_{0}$ for the ordered pair
$(W_{1},W_{2})$. 
\et

{\bf Proof}. Let $W$ be a $V$-module and let $I(\cdot,x)$ be any
intertwining operator of type
$\left(\!\begin{array}{c}W\\W_{1} W_{2}\end{array}\!\right)$. Let $D_{i}$ 
be the projection 
of $W$ onto $W^{(i)}$ for $i=1,\cdots,k$. 
Then $D_{i}I(\cdot,x)$ is an intertwining operator of type
$\left(\!\begin{array}{c}W^{(i)}\\W_{1} W_{2}\end{array}\!\right)$. 
By Corollary \ref{c3.16}, we obtain a $V$-homomorphism $g_{i}$ from  
$T(W_{1},W_{2})$ to $W^{(i)}$ satisfying the condition: 
\begin{eqnarray}
g_{i}\pi (Y_{t}(u_{1},x)\otimes u_{2})=D_{i}I^{o}(u_{1},x)u_{2}\;\;\;\;\mbox{for }
u_{1}\in W_{1}, u_{2}\in W_{2}.
\end{eqnarray}
Since $g_{i}P_{j}=0$ for $j\ne i$, we obtain
$g_{i}\circ F(u_{1},x)u_{2}=D_{i}I(u_{1},x)u_{2}$ for $u_{1}\in W_{1},
u_{2}\in W_{2}$.  Set $g=g_{1}\oplus \cdots \oplus g_{k}$. Then
$$g\circ F(u_{1},x)u_{2}=I(u_{1},x)u_{2}\;\;\;\mbox{ for }u_{1}\in W_{1},
u_{2}\in W_{2}.$$
{}From the construction of $T(W_{1},W_{2})$, $F(\cdot,x)$ is surjective 
in the sense of Lemma \ref{l3.3}, {i,e.,} all the coefficients of 
$F(u_{1},x)u_{2}$
for $u_{i}\in W_{i}$ linearly span $T(W_{1},W_{2})$. Thus such a $g$ is unique.
Then the pair $(T(W_{1},W_{2}), F(\cdot,x))$ 
is a tensor product for the ordered pair
$(W_{1},W_{2})$. $\;\;\;\;\Box$


\section{An analogue of the ``Hom''-functor and a generalized nuclear
democracy theorem}
In this section we shall introduce the notion of what we call ``generalized
intertwining operator'' from a $V$-module $W_{1}$ to another
$V$-module $W_{2}$. The notion of generalized intertwining operator can be
considered as a generalization of the physicists' notion 
of ``primary field''
(cf. [BPZ], [MS] and [TK]) to the notion of general (non-primary) field.
On the other hand, it exactly reflects the main features of $I(u,x)$ for
$u\in M$, where $M$ is a $V$-module and $I(\cdot,x)$ is an
intertwining operator of type
$\left(\!\begin{array}{c}W_{2}\\MW_{1}\end{array}\!\right)$. We prove that 
$G(W_{1},W_{2})$, the space of all generalized intertwining operators,
is a weak $V$-module (Theorem \ref{t4.6}), which satisfies a
certain universal property in terms of 
the space of intertwining operators of a certain type (Theorem \ref{t4.8}).
If the vertex operator algebra $V$ satisfies certain finiteness
and semisimplicity conditions, we prove that there exists a unique
maximal submodule 
$\Delta(W_{1},W_{2})$ of $G(W_{1},W_{2})$ so that the contragredient
module of $\Delta(W_{1},W_{2})$ is a tensor product module for
the ordered pair $(W_{1}, W_{2}')$ (Theorem \ref{t4.13}). Using 
Theorem \ref{t4.8} we derive a generalized form of 
the nuclear democracy theorem of Tsuchiya and Kanie [TK] (Theorem \ref{t4.14}).
All these results show that the notion of $G(W_{1},W_{2})$ is 
an analogue of the classical ``${\rm Hom}$''-functor.

{\em Throughout this section, $V$ will be a fixed vertex operator algebra.}

\bd{d4.1} Let $W_{1}$ and $W_{2}$ be $V$-modules.  A
{\em generalized intertwining operator} from $W_{1}$ 
to $W_{2}$ is an element 
$\Phi (x)=\displaystyle{\sum_{\alpha\in {\C}}\Phi_{\alpha}x^{-\alpha-1}}
\in \left({\rm Hom}(W_{1},W_{2})\right)\{x\}$ satisfying the following
conditions (G1)-(G3): 

(G1)\hspace{0.25cm}  For any $\alpha\in {\C}, u_{1}\in W_{1}$,
$\Phi_{\alpha+n}u_{1}=0$ for $n\in {\Z}$ sufficiently large;

(G2)\hspace{0.25cm} $[L(-1),\Phi (x)]=\Phi'(x)
 \left(=\displaystyle{{d\over dx}}\Phi(x)\right);$

(G3)\hspace{0.25cm} For any $a\in V$, there exists a positive integer
$k$ such that 
\begin{eqnarray}\label{el}
(x_{1}-x_{2})^{k}Y_{W_{2}}(a,x_{1})\Phi(x_{2})
=(x_{1}-x_{2})^{k}\Phi(x_{2})Y_{W_{1}}(a,x_{1}).
\end{eqnarray}
\ed
Denote by ${\rm G}(W_{1},W_{2})$ the space of all
generalized intertwining operators from $W_{1}$ to
$W_{2}$.
A generalized intertwining  operator $\Phi(x)$ is said to be
{\it homogeneous of weight $h$} if it satisfies the following condition:
\begin{eqnarray}\label{eh}
[L(0),\Phi(x)]=\left(h+x{d\over dx}\right)\Phi(x). 
\end{eqnarray}
A generalized intertwining operator $\Phi(x)$ of weight $h$ is said to 
be {\em primary} if the following condition holds:
\begin{eqnarray}\label{ep2}
[L(m), \Phi(x)]=x^{m}\left((m+1)h+x{d\over dx}\right)\Phi(x)\;\;\;\mbox{ for }m\in {\Z}.
\end{eqnarray}
Denote by ${\rm G}(W_{1},W_{2})_{(h)}$ the space of all weight-$h$
homogeneous generalized intertwining operators from $W_{1}$ to
$W_{2}$ and set
\begin{eqnarray}
{\rm G}^{0}(W_{1},W_{2})=\oplus_{h\in {\C}}{\rm G}(W_{1},W_{2})_{(h)}.
\end{eqnarray}

Let $W(W_{1},W_{2})$ be the space consisting of each element
$\Phi(x)\in \left({\rm Hom}_{{\C}}(W_{1},W_{2})\right)$
which satisfies the condition (G1) and let $E(W_{1},W_{2})$
be the space consisting of each element
$\Phi(x)\in \left({\rm Hom}_{{\C}}(W_{1},W_{2})\right)$
which satisfies the conditions (G1) and (G3).
For any $a\in V$, we
define the left and the right actions of $\hat{V}$ on  $W(W_{1},W_{2})$
as follows:
\begin{eqnarray}
Y_{t}(a,x_{0})*\Phi(x_{2}):&=&{\rm Res}_{x_{1}}x_{0}^{-1}\delta\left(
\frac{x_{1}-x_{2}}{x_{0}}\right)Y_{W_{2}}(a,x_{1})\Phi(x_{2})\\
&=&Y_{W_{2}}(a,x_{0}+x_{2})\Phi(x_{2}).\\
\Phi(x_{2})* Y_{t}(a,x_{0}):&=&{\rm Res}_{x_{1}}x_{0}^{-1}\delta\left(
\frac{-x_{2}+x_{1}}{x_{0}}\right)\Phi(x_{2})Y_{W_{1}}(a,x_{1})\\
&=&\Phi(x_{2})(Y_{W_{1}}(a,x_{0}+x_{2})-Y_{W_{1}}(a,x_{2}+x_{0})).
\end{eqnarray}

\bp{p4.2}  a) $W(W_{1},W_{2})$ is a left $
g(V)$-module of level one under the defined left action.

b) $W(W_{1},W_{2})$ is a right $g(V)$-module of level zero under
the defined right action.
\ep

{\bf Proof.} a) First we check that $W(W_{1},W_{2})$ is closed under
the left action. For any $a\in V, m\in {\Z}, \Phi(x)\in
W(W_{1},W_{2}), u\in W_{1}$, by definition we have:
\begin{eqnarray}
\left((t^{m}\otimes a)*\Phi(x)\right)u&=&{\rm Res}_{x_{0}}x_{0}^{m}
Y_{W_{2}}(a,x_{0}+x)\Phi(x_{2})u\nonumber\\
&=&\sum_{i=0}^{\infty}
{-m+i-1\choose i}x^{i}a_{m-i}\Phi(x)u.
\end{eqnarray}
Then it is clear that $(t^{m}\otimes a)*\Phi(x_{2})$ satisfies (G1). 
Next, we check the defining relations for $g(V)$. By
definition we have
\begin{eqnarray}
Y_{t}({\bf 1},x_{0})*\Phi(x_{2})=Y_{W_{2}}({\bf 1},x_{0}+x_{2})\Phi(x_{2})
=\Phi(x_{2})
\end{eqnarray}
and
\begin{eqnarray}
Y_{t}(L(-1)a,x_{0})*\Phi(x_{2})&=&Y_{W_{2}}(L(-1)a,x_{0}+x_{2})\Phi(x_{2})
\nonumber\\
&=&{\partial\over \partial x_{0}}Y_{W_{2}}(a,x_{0}+x_{2})\Phi(x_{2})
\nonumber\\
&=&{\partial\over \partial x_{0}}Y_{t}(a,x_{0})*\Phi(x_{2}).
\end{eqnarray}
Furthermore, for any $a,b\in V$, we have
\begin{eqnarray}
& &Y_{t}(a,x_{1})*Y_{t}(b,x_{2})*\Phi(x_{3})\nonumber\\
&=&Y_{t}(a,x_{1})*\left(Y_{W_{2}}(b,x_{2}+x_{3})\Phi(x_{3})\right)
\nonumber\\
&=&Y_{W_{2}}(a,x_{1}+x_{3})Y_{W_{2}}(b,x_{2}+x_{3})\Phi(x_{3}).
\end{eqnarray}
Similarly, we have
\begin{eqnarray}
Y_{t}(b,x_{2})*Y_{t}(a,x_{1})*\Phi(x_{3})=Y_{W_{2}}(b,x_{2}+x_{3})
Y_{W_{2}}(a,x_{1}+x_{3})\Phi(x_{3}).
\end{eqnarray}
Therefore
\begin{eqnarray}
& &Y_{t}(a,x_{1})*Y_{t}(b,x_{2})*\Phi(x_{3})-Y_{t}(b,x_{2})*Y_{t}(a,x_{1})*
\Phi(x_{3})\nonumber\\
&=&{\rm Res}_{x_{0}}(x_{2}+x_{3})^{-1}\delta\left(\frac{x_{1}+x_{3}-x_{0}}
{x_{2}+x_{3}}\right)Y_{W_{2}}(Y(a,x_{0})b,x_{2}+x_{3})\Phi(x_{3})
\nonumber\\
&=&{\rm Res}_{x_{0}}x_{1}^{-1}\delta\left(\frac{x_{2}+x_{0}}{x_{1}}\right)
Y_{W_{2}}(Y(a,x_{0})b,x_{2}+x_{3})\Phi(x_{3})\nonumber\\
&=&{\rm Res}_{x_{0}}x_{2}^{-1}\delta\left(\frac{x_{1}-x_{0}}{x_{2}}\right)
Y_{t}(Y(a,x_{0})b,x_{2})*\Phi(x_{3}).
\end{eqnarray}
Then a) is  proved.

The proof of b) is similar to the proof of a), but for completeness, we also 
write the details. For any $a\in V, \Phi(x)\in W(W_{1},W_{2})$, by
definition we have
\begin{eqnarray}
& &\Phi(x_{2})*Y_{t}(L(-1)a,x_{0})\nonumber\\
&=&\Phi(x_{2})(Y_{W_{1}}(L(-1)a,x_{0}+x_{2})-Y_{W_{1}}(L(-1)a,x_{2}+x_{0}))
\nonumber\\
&=&{\partial\over\partial x_{0}}\left(\Phi(x_{2})(Y_{W_{1}}(a,x_{0}+x_{2})
-Y_{W_{1}}(a,x_{2}+x_{0})\right)\nonumber\\
&=&{\partial\over\partial x_{0}}\Phi(x_{2})*Y_{t}(a,x_{0}).
\end{eqnarray}
For any $a,b\in V$, we have
\begin{eqnarray}
& &\Phi(x_{3})*Y_{t}(a,x_{1})*Y_{t}(b,x_{2})\nonumber\\
&=&{\rm Res}_{x_{4}}x_{1}^{-1}\delta\left(\frac{-x_{3}+x_{4}}{x_{1}}\right)
\left(\Phi(x_{3})Y_{W_{1}}(a,x_{4})\right)*Y_{t}(b,x_{2})\nonumber\\
&=&{\rm Res}_{x_{4}}{\rm Res}_{x_{5}}x_{1}^{-1}\delta\left(\frac{-x_{3}+x_{4}}
{x_{1}}\right)x_{2}^{-1}\delta\left(\frac{-x_{3}+x_{5}}
{x_{2}}\right)\Phi(x_{3})Y_{W_{1}}(a,x_{4})Y_{W_{1}}(b,x_{5}).\nonumber\\
& &\mbox{}
\end{eqnarray}
Similarly, we have
\begin{eqnarray}
& &\Phi(x_{3})*Y_{t}(b,x_{2})*Y_{t}(a,x_{1})\nonumber\\
&=&{\rm Res}_{x_{4}}{\rm Res}_{x_{5}}x_{1}^{-1}\delta\left(\frac{-x_{3}+x_{4}}
{x_{1}}\right)x_{2}^{-1}\delta\left(\frac{-x_{3}+x_{5}}
{x_{2}}\right)\Phi(x_{3})Y_{W_{1}}(b,x_{5})Y_{W_{1}}(a,x_{4}).\nonumber\\
& &\mbox{}
\end{eqnarray}
Thus 
\begin{eqnarray}
& &\Phi(x_{3})*Y_{t}(a,x_{1})*Y_{t}(b,x_{2})-\Phi(x_{3})*Y_{t}(b,x_{2})*
Y_{t}(a,x_{1})\nonumber\\
&=&{\rm Res}_{x_{4}}{\rm Res}_{x_{5}}{\rm Res}_{x_{0}}x_{1}^{-1}\delta\left(
\frac{-x_{3}+x_{4}}{x_{1}}\right)x_{2}^{-1}\delta\left(\frac{-x_{3}+x_{5}}
{x_{2}}\right)x_{5}^{-1}\delta\left(\frac{x_{4}-x_{0}}{x_{5}}\right)\nonumber\\
& &\cdot \Phi(x_{3})Y_{W_{1}}(Y(a,x_{0})b,x_{5})\nonumber\\
&=&{\rm Res}_{x_{5}}{\rm Res}_{x_{0}}x_{1}^{-1}\delta\left(
\frac{-x_{3}+x_{5}+x_{0}}{x_{1}}\right)x_{2}^{-1}\delta\left(
\frac{-x_{3}+x_{5}}{x_{2}}\right)\Phi(x_{3})Y_{W_{1}}(Y(a,x_{0})b,x_{5})
\nonumber\\
&=&{\rm Res}_{x_{5}}{\rm Res}_{x_{0}}x_{1}^{-1}\delta\left(
\frac{x_{2}+x_{0}}{x_{1}}\right)x_{2}^{-1}\delta\left(
\frac{-x_{3}+x_{5}}{x_{2}}\right)\Phi(x_{3})Y_{W_{1}}(Y(a,x_{0})b,x_{5})
\nonumber\\
&=&\Phi(x_{3})*{\rm Res}_{x_{0}}x_{2}^{-1}\delta\left(\frac{x_{1}-x_{0}}
{x_{2}}\right)Y_{t}(Y(a,x_{0})b,x_{2}).
\end{eqnarray}
Then the proof is complete.$\;\;\;\;\Box$

For any $a\in V, \Phi(x)\in W(W_{1},W_{2})$, we define               
\begin{eqnarray}
& &\!\!Y_{t}(a,x_{0})\circ \Phi(x_{2}):=\nonumber\\
&=&\!\!Y_{t}(a,x_{0})*\Phi(x_{2})-\Phi(x_{2})*Y_{t}(a,x_{0})\nonumber\\
&=&\!\!{\rm Res}_{x_{1}}\!\left(\!x_{0}^{-1}\delta\!\left(\!\frac{x_{1}-x_{2}}{x_{0}}\!
\right)Y_{W_{2}}(a,x_{1})\Phi (x_{2})-
x_{0}^{-1}\delta\!\left(\!\frac{x_{2}-x_{1}}{-x_{0}}\!\right)\Phi (x_{2})
Y_{W_{1}}(a,x_{1})\!\right)\!,
\end{eqnarray}
or equivalently
\begin{eqnarray}\label{e4.21}
a(m)\circ \Phi (x_{2})={\rm Res}_{x_{1}}\!\left(\! (x_{1}-x_{2})^{m}
Y_{W_{2}}(a,x_{1})\Phi (x_{2})-(-x_{2}+x_{1})^{m}\Phi(x_{2})Y_{W_{1}}(a,x_{1})
\!\right)
\end{eqnarray}
for any $m\in {\Z}$. From the classical Lie algebra theory, we have

\bc{c4.3} Under the defined action $''\circ''$,
$W(W_{1},W_{2})$ becomes a $g(V)$-module (of level one).
\ec

\bl{l4.4}
Let $\Phi (x)\in W(W_{1},W_{2})$
satisfying (\ref{eh}) for some complex number $h$ and let $a$ be any
homogeneous element of $V$. Then
\begin{eqnarray}
[L(0), Y_{t}(a,x_{0})\circ \Phi (x_{2})]
=\left({\rm wt}a+h+x_{0}{\partial\over\partial
x_{0}}+x_{2}{\partial\over\partial x_{2}}\right)Y_{t}(a,x_{0})\circ
\Phi (x_{2}).
\end{eqnarray}
\el

{\bf Proof.} By definition we have
\begin{eqnarray}
& &[L(0),Y_{t}(a,x_{0})\circ \Phi (x_{2})]\nonumber\\
&=&{\rm Res}_{x_{1}}x_{0}^{-1}\delta\left(\frac{x_{1}-x_{2}}{x_{0}}\right)
[L(0),Y(a,x_{1})\Phi (x_{2})]\nonumber\\
& &-{\rm Res}_{x_{1}}x_{0}^{-1}\delta\left(\frac{-x_{2}+x_{1}}{x_{0}}\right)
[L(0),\Phi (x_{2})Y(a,x_{1})]\nonumber\\
&=&{\rm Res}_{x_{1}}x_{0}^{-1}\delta\left(\frac{x_{1}-x_{2}}{x_{0}}\right)
\left({\rm wt}a+x_{1}{\partial\over\partial
x_{1}}+h+x_{2}{\partial\over\partial x_{2}}\right)Y(a,x_{1})\Phi (x_{2})
\nonumber\\
& &-{\rm Res}_{x_{1}}x_{0}^{-1}\delta\left(\frac{x_{2}-x_{1}}{-x_{0}}\right)
\left({\rm wt}a+x_{1}{\partial\over\partial
x_{1}}+h+x_{2}{\partial\over\partial x_{2}}\right)\Phi (x_{2})Y(a,x_{1})
\nonumber\\
&=&({\rm wt}a+h)Y_{t}(a,x_{0})\circ \Phi (x_{2})\nonumber\\
& &-{\rm Res}_{x_{1}}\left({\partial\over\partial x_{1}}x_{1}x_{0}^{-1}
\delta\left(\frac{x_{1}-x_{2}}{x_{0}}\right)\right)Y(a,x_{1})\Phi (x_{2})
\nonumber\\
& &+{\rm Res}_{x_{1}}x_{0}^{-1}\delta\left(\frac{x_{1}-x_{2}}{x_{0}}\right)
x_{2}{\partial\over\partial x_{2}}Y(a,x_{1})\Phi (x_{2})\nonumber\\
& &-{\rm Res}_{x_{1}}x_{0}^{-1}\delta\left(\frac{x_{2}-x_{1}}{-x_{0}}\right)
x_{2}{\partial\over\partial x_{2}}\Phi (x_{2})Y(a,x_{1})\nonumber\\
& &+{\rm Res}_{x_{1}}\left({\partial\over\partial x_{1}}x_{1}x_{0}^{-1}\delta
\left(\frac{x_{2}-x_{1}}{-x_{0}}\right)\right)\Phi (x_{2})Y(a,x_{1}).
\end{eqnarray}
Since
\begin{eqnarray}
{\partial\over\partial x_{0}}x_{0}^{-1}\delta\left(
\frac{x_{1}-x_{2}}{x_{0}}\right)
={\partial\over\partial x_{2}}x_{0}^{-1}\delta\left(
\frac{x_{1}-x_{2}}{x_{0}}\right)
=-{\partial\over\partial x_{1}}x_{0}^{-1}\delta\left(
\frac{x_{1}-x_{2}}{x_{0}}\right),
\end{eqnarray}
we have
\begin{eqnarray}
& &{\partial\over\partial x_{1}}\left(x_{1}x_{0}^{-1}\delta\left(
\frac{x_{1}-x_{2}}{x_{0}}\right)\right)\nonumber\\
&=&{\partial\over\partial x_{1}}\left((x_{0}+x_{2})x_{0}^{-1}\delta\left(
\frac{x_{1}-x_{2}}{x_{0}}\right)\right)\nonumber\\
&=&x_{0}{\partial\over\partial x_{1}}x_{0}^{-1}\delta\left(
\frac{x_{1}-x_{2}}{x_{0}}\right)+x_{2}{\partial\over\partial x_{1}}x_{0}^{-1}
\delta\left(\frac{x_{1}-x_{2}}{x_{0}}\right)\nonumber\\
&=&-x_{0}{\partial\over\partial x_{0}}x_{0}^{-1}\delta\left(
\frac{x_{1}-x_{2}}{x_{0}}\right)
-x_{2}{\partial\over\partial x_{2}}x_{0}^{-1}\delta\left(
\frac{x_{1}-x_{2}}{x_{0}}\right).
\end{eqnarray}
Similarly, we have
\begin{eqnarray}
{\partial\over\partial x_{1}}\left(x_{1}x_{0}^{-1}\delta\left(
\frac{x_{2}-x_{1}}{-x_{0}}\right)\right)
=-\left(x_{0}{\partial\over\partial x_{0}}+x_{2}{\partial\over\partial x_{2}}
\right)x_{1}x_{0}^{-1}\delta\left(\frac{x_{2}-x_{1}}{-x_{0}}\right).
\end{eqnarray}
Therefore, we obtain
\begin{eqnarray}
& &[L(0),Y_{t}(a,x_{0})\circ \Phi(x_{2})]\nonumber\\
&=&\left({\rm wt}a+h+x_{0}{\partial\over\partial x_{0}}+
x_{2}{\partial\over\partial x_{2}}\right)Y_{t}(a,x_{0})\circ \Phi(x_{2}).
\;\;\;\;\Box \end{eqnarray}

\bp{p4.5} The subspaces $\!E(\!W_{1},W_{2}\!)\!$ and $\!{\rm G}(\!W_{1},W_{2}\!)$ are
restricted $g(V)$-submodules of $W(W_{1},W_{2})$ and ${\rm
G}^{0}(W_{1},W_{2})$  is a ${\C}$-graded $g(V)$-module.
\ep

{\bf Proof.} For any $a\in V, m\in {\Z},\Phi(x)\in E(W_{1},W_{2})$,
it follows from  the proof of Proposition 3.2.7 in [L1] (for 
an analogous result) that $a(m)\circ \Phi(x)\in E(W_{1},W_{2})$. Thus $E(W_{1},W_{2})$
is a submodule of $W(W_{1},W_{2})$.
For $\Phi(x)\in G(W_{1},W_{2})$, since
\begin{eqnarray}
& &[L(-1), Y_{t}(a,x_{0})*\Phi(x_{2})]\nonumber\\
&=&[L(-1),Y_{W_{2}}(a,x_{0}+x_{2})\Phi(x_{2})]\nonumber\\
&=&[L(-1),Y_{W_{2}}(a,x_{0}+x_{2})]\Phi(x_{2})+Y_{W_{2}}(a,x_{0}+x_{2})
[L(-1),\Phi(x_{2})]\nonumber\\
&=&{\partial\over \partial x_{2}}\left(Y_{W_{2}}(a,x_{0}+x_{2})\Phi(x_{2})
\right)\nonumber\\
&=&{\partial\over \partial x_{2}}Y_{t}(a,x_{0})*\Phi(x_{2}),
\end{eqnarray}
$a(m)*\Phi(x_{2})$ satisfies (G2). Thus $a(m)*\Phi(x_{2})\in G(W_{1},W_{2})$.
Similarly, since
\begin{eqnarray}
& &{\partial\over \partial x_{2}}\left(\Phi(x_{2})*Y_{t}(a,x_{0})\right)
\nonumber\\
&=&{\partial\over \partial x_{2}}\left(\Phi(x_{2})(Y_{W_{1}}(a,x_{0}+x_{2})
-Y_{W_{1}}(a,x_{2}+x_{0}))\right)\nonumber\\
&=&\Phi'(x_{2})(Y_{W_{1}}(a,x_{0}+x_{2})
-Y_{W_{1}}(a,x_{2}+x_{0}))\nonumber\\
& &+\Phi(x_{2})(Y_{W_{1}}(L(-1)a,x_{0}+x_{2})
-Y_{W_{1}}(L(-1)a,x_{2}+x_{0}))\nonumber\\
&=&[L(-1),\Phi(x_{2})*Y_{t}(a,x_{0})],
\end{eqnarray}
we obtain $\Phi(x_{2})* a(m)\in G(W_{1},W_{2})$. Therefore $a(m)\circ \Phi(x)\in G(W_{1},W_{2})$.
Thus ${\rm G}(W_{1},W_{2})$ is a submodule. That ${\rm
G}^{0}(W_{1},W_{2})$  is a ${\C}$-graded $g(V)$-module follows from Lemma \ref{l4.4}. It follows from (\ref{e4.21}) and (G3) that $E(W_{1},W_{2})$ is a restricted $g(V)$-module and so are $G(W_{1},W_{2})$ and $G^{0}(W_{1},W_{2})$.
Then the proof is complete.
$\;\;\;\;\Box$

Define a linear map $F(\cdot,x)$ from $E(W_{1},W_{2})$ to ${\rm Hom}(W_{1},W_{2})\{x\}$ as follows:
\begin{eqnarray}
F(\Phi,x)u_{1}=\Phi (x)u_{1}\;\;\;\mbox{ for }\Phi \in E(W_{1},W_{2}), u_{1}\in W_{1}.
\end{eqnarray}
For $a\in V, \Phi \in E(W_{1},W_{2})$, we have
\begin{eqnarray}
& &F(Y(a,x_{0})\Phi,x_{2})\nonumber\\
&=&(Y(a,x_{0})\Phi)(x)|_{x=x_{2}}\nonumber\\
&=&{\rm Res}_{x_{1}}\left(x_{0}^{-1}\delta\left(\frac{x_{1}-x}{x_{0}}
\right)Y(a,x_{1})\Phi(x)-x_{0}^{-1}\delta\left(\frac{-x+x_{1}}
{x_{0}}\right)\Phi(x)Y(a,x_{1})\right)|_{x=x_{2}}\nonumber\\
&=&{\rm Res}_{x_{1}}\left(x_{0}^{-1}\delta\!\left(\frac{x_{1}-x_{2}}{x_{0}}
\right)Y(a,x_{1})F(\Phi,x_{2})-x_{0}^{-1}\delta\!\left(\frac{x_{2}-x_{1}}
{-x_{0}}\right)F(\Phi,x_{2})Y(a,x_{1})\right).\nonumber\\
& &\mbox{}
\end{eqnarray}
It is well-known that this iterate formula implies the associativity [FHL]. Furthermore, 
(G3) gives the commutativity for $F(\cdot,x)$. Therefore,
$F(\cdot,x)$ satisfies the Jacobi 
identity ([DL], [FHL], [L1]).
Thus $F(\cdot,x)$ is a weak intertwining operator.
It is clear that $F(\cdot,x)$ is injective in the sense that $F(\Phi,x)=0$ implies
$\Phi=0$ for $\Phi\in E(W_{1},W_{2})$.
Furthermore, if $\Phi(x)\in G(W_{1},W_{2})$, by definition we have
\begin{eqnarray}
& &F(L(-1)\Phi,x)u_{1}=(L(-1)\Phi )(x)u_{1}=(L(-1)\circ \Phi(x))u_{1}={d\over
dx}\Phi(x)u_{1}={d\over dx}F(\Phi,x)u_{1}.\nonumber\\ 
& &\mbox{}
\end{eqnarray}
Therefore, $F(\cdot,x)$  is an intertwining 
operator of
type $\left(\!\begin{array}{c}W_{2}\\G(W_{1},W_{2}) W_{1}\end{array}\!\right)$ after 
restricted to $G(W_{1},W_{2})$.

\bt{t4.6} The $g(V)$-module $E(W_{1},W_{2})$ and $G(W_{1},W_{2})$ are
weak $V$-modules.
\et

{\bf Proof }\footnote{This was proved directly in [L1].}. By Proposition \ref{p3.11}, we get
$$F(\Phi,x)=0\;\;\;\mbox{ for any }\Phi \in J(E(W_{1},W_{2})).$$
Since $F(\cdot,x)$ injective, $J(E(W_{1},W_{2}))=0$. That is, $E(W_{1},W_{2})$ is a 
weak $V$-module.$\;\;\;\;\Box$

Let $M$ be another $V$-module and let $f\in {\rm Hom}_{V}(M,
G(W_{1},W_{2}))$. Then $F(f\cdot,x)$ is an intertwining operator of type
$\left(\!\begin{array}{c}W_{2}\\M W_{1}\end{array}\!\right)$. 
Since $F(\cdot,x)$ is injective, we obtain an injective linear map
\begin{eqnarray}
\theta : & &{\rm Hom}_{V}(M,G(W_{1},W_{2}))\rightarrow I\!\left(\!
\begin{array}{c}W_{2}\\M W_{1}\end{array}\!\right)\nonumber\\
& &f\mapsto F(f\cdot,x).
\end{eqnarray}

On the other hand, for any intertwining operator $I(\cdot,x)$ of type
 $\left(\!\begin{array}{c}W_{2}\\M W_{1}\end{array}\!\right)$, it is clear
that $I(u,x)\in G(W_{1},W_{2})$ for any $u\in M$. Then we
obtain a linear map $f_{I}$ from $M$ to $G(W_{1},W_{2})$
defined by $f_{I}(u)=I(u,x)$. For any $a\in V, u\in M$, from the definition of 
$Y(a,x_{0})\circ I(u,x)$ we get
\begin{eqnarray}
& &f_{I}(Y(a,x_{0})u)\nonumber\\
&=&I(Y(a,x_{0})u,x)\nonumber\\
&=&{\rm Res}_{x_{1}}\left(x_{0}^{-1}\delta\!\left(\frac{x_{1}-x}{x_{0}}
\right)Y(a,x_{1})I(u,x)-x_{0}^{-1}\delta\!\left(\frac{x-x_{1}}
{-x_{0}}\right)I(u,x)Y(a,x_{1})\right)\nonumber\\
&=&Y(a,x_{0})\circ I(u,x)\nonumber\\
&=&Y(a,x_{0})f_{I}(u).
\end{eqnarray}
Thus $f_{I}$ is a $V$-homomorphism such that $F(f_{I}\cdot,x)=I(\cdot,x)$.
Since $F(\cdot,x)$ is injective, such an $f_{I}$ is unique. Therefore we obtain

\bt{t4.8} 
Let $W_{1}$and $W_{2}$ be $V$-modules. Then
(a) For any $V$-module $M$ and any intertwining operator $I(\cdot,x)$
of type $\left(\!\begin{array}{c}W_{2}\\M W_{1}\end{array}\!\right)$,
there exists a unique $V$-homomorphism $f$ from $M$ to 
$G(W_{1},W_{2})$ such that $I(u,x)=F(f(u),x)$ for $u\in M$.

(b) The linear space ${\rm Hom}_{V}(M, G(W_{1},W_{2}))$ is naturally
isomorphic to
$I\!\left(\!\begin{array}{c}W_{2}\\M W_{1} \end{array}\!\right)$ for any 
$V$-module $M$.
\et

The universal property in Theorem \ref{t4.8} looks very
much like the universal property for a tensor product in Definition
\ref{d3.1} and also in [HL1]. 
Next, we study the relation between
$G(W_{1},W_{2})$ and the contragredient module
of tensor product of $W_{1}$ and $W_{2}'$. 

\br{r4.11}Let $M$ be any $V$-module. Then it was proved in [L1] that
$G(V,M)\simeq M$. If $M=V$, then $V=G(V,V)$. That is, any
generalized intertwining operator is a vertex operator. In this
special case, this has been proved in [G].
\er

For any two $V$-modules $W_{1}$ and $W_{2}$, let $\Delta
(W_{1},W_{2})$ be the sum of all $V$-modules inside the 
weak $V$-module $G(W_{1},W_{2})$.

\bp{p4.12}\footnote{A similar result has also been obtained in [HL0-4].}
Let $V$ be a vertex operator
algebra satisfying the following conditions: (1) There are finitely
many inequivalent irreducible $V$-modules. (2) Any $V$-module is
completely reducible. (3)  Any fusion rule for three modules is 
finite.  Then for any $V$-modules $W_{1}$ and $W_{2}$,
$\Delta(W_{1},W_{2})$ is the unique maximal $V$-module inside
the weak module $G(W_{1},W_{2})$. 
\ep

{\bf Proof.} It follows from the condition (2) that
$\Delta(W_{1},W_{2})$ is a direct sum of irreducible $V$-modules. It
follows from Theorem \ref{t4.8} and the condition (3) that the multiplicity
of each irreducible $V$-module in  $\Delta(W_{1},W_{2})$ is finite.
Therefore $\Delta(W_{1},W_{2})$ is a direct sum of finitely many
irreducible $V$-modules. That is, $\Delta(W_{1},W_{2})$ is a
$V$-module. By the definition of $\Delta(W_{1},W_{2})$, it is clear
that $\Delta(W_{1},W_{2})$ is the unique 
maximal $V$-module inside the weak $V$-module ${\rm G}(W_{1},W_{2})$. 
$\;\;\;\;\Box$

Let $V$ be a vertex operator algebra satisfying the conditions (1)-(3)
of Proposition \ref{p4.12}  and let $W_{1}$ and $W_{2}$ be any two
$V$-modules. Let $\bar{F}(\cdot,x)$ be
the restriction of $F(\cdot,x)$ on $\Delta(W_{1},W_{2})$ so that $\bar{F}(\cdot,x)$ is 
an intertwining operator of type 
$\left(\!\begin{array}{c}W_{2}\\
\Delta(W_{1},W_{2}) W_{1}\end{array}\!\right)$ such that
\begin{eqnarray}
\bar{F}(\Phi ,x)=\Phi (x)\;\;\;\mbox{ for any }\Phi \in \Delta (W_{1},W_{2}).
\end{eqnarray}
Then by Proposition \ref{p2.7},  the transpose operator $\bar{F}^{t}(\cdot,x)$ of
$\bar{F}(\cdot,x)$ is an
intertwining operator of type $\left(\!\begin{array}{c}W_{2}\\
W_{1}\Delta(W_{1},W_{2})\end{array}\!\right)$. Furthermore, it follows
{}from Proposition \ref{p2.7} that $(\bar{F}^{t})'(\cdot,x)$ is an intertwining
operator of type $\left(\!\begin{array}{c}(\Delta(W_{1},W_{2}))'\\ W_{1}
W_{2}'\end{array}\!\right)$.

\bt{t4.13} If $V$ satisfies the conditions (1)-(3) of 
Proposition \ref{p4.12}, then the pair
$((\Delta (W_{1},W_{2})', (\bar{F}^{t})'(\cdot,x))$ is a tensor product
for the ordered pair $(W_{1},W_{2}')$ in the category of $V$-modules.
\et

{\bf Proof.} Let $W$ be any $V$-module and let $I(\cdot,x)$ be any
intertwining operator of type $\left(
\!\begin{array}{c}W\\W_{1}W_{2}'\end{array}\!\right)$. It follows from 
Proposition \ref{p2.7} that
 $(I')^{t}(\cdot,x)$ is an intertwining operator
of type $\left(\!
\begin{array}{c}W_{2}\\W'W_{1}\end{array}\!\right)$. From Theorem \ref{t4.6},
there exists a (unique) $V$-homomorphism $\psi$ from $W'$ to $
G(W_{1},W_{2})$ such that $(I')^{t}(w',x)=\bar{F}(\psi (w'),x)$ for any $w'\in M'$.
It follows from the definition of $\Delta(W_{1},W_{2})$ that $\psi$ is
a $V$-homomorphism from $W'$ to $\Delta(W_{1},W_{2})$. Therefore,
we obtain a $V$-homomorphism $\psi '$ from $(\Delta(W_{1},W_{2}))'$ to $W$.
For any $w'\in W', u_{1}\in W_{1}, u_{2}'\in W_{2}'$, by using FLM's 
conjugation formulas [FHL] we obtain
\begin{eqnarray}
& &\<w', \psi' (\bar{F}^{t})'(u_{1},x)u_{2}'\>\nonumber\\
&=&\<\bar{F}^{t}(e^{xL(1)}(e^{\pi i}x^{-2})^{L(0)}u_{1},x^{-1})\psi w',u_{2}'\>
\nonumber\\
&=&\<I'(e^{xL(1)}(e^{\pi i}x^{-2})^{L(0)}u_{1},x^{-1}) w',u_{2}'\>\nonumber\\
&=&\<w',I(e^{x^{-1}L(1)}(e^{\pi i}x^{2})^{L(0)}e^{xL(1)}
(e^{\pi i}x^{-2})^{L(0)}u_{1},x)u_{2}'\>\nonumber\\
&=&\<w', I(e^{2\pi i L(0)}u_{1},x)u_{2}'\>.
\end{eqnarray}
For any $V$-module $M$, we define a linear endomorphism $t_{M}$ of $M$ by:
$t_{M}(u)=e^{2\pi iL(0)}u$ for $u\in M$. Then one can easily prove that 
$t_{M}$ is a $V$-automorphism of $M$ so that $t_{M}$ is a scalar if $M$ is 
irreducible. Let $t_{W_{1}}=\alpha$. Then 
$\alpha^{-1}\psi'(\bar{F}^{t})'(\cdot,x)=I(\cdot,x)$.
The uniqueness of $\alpha^{-1}\psi'$ follows from the uniqueness of $\psi$.
Then the proof is complete.$\;\;\;\;\Box$

\begin{rem}
It was was proved in [DLM2] that the category $\cal{C}$ of all weak 
$V$-modules is a semisimple category
for vertex operator algebras
$L(\ell,0)$, associated to an integrable highest weight module of level 
$\ell$ for an affine Lie algebra, $L(c_{p,q},0)$, associated to the 
irreducible highest weight module for the Virasoro algebra with central 
charge $c_{p,q}=1-\frac{6(p-q)^{2}}{pq}$, 
the moonshine module vertex operator algebra $V^{\natural}$ and $V_{L}$,
associated to any even positive-definite lattice $L$.
Thus, for these vertex operator algebras, we have
$G(W_{1},W_{2})=\Delta(W_{1},W_{2})$. 
\end{rem}

Let $U$ be an irreducible $g(V)_{0}$-module on which $L(0)$ acts as a 
scalar $h$.
Define $g(V)_{-}U=0$. Then $U$ becomes a $(g(V)_{-}+g(V)_{0})$-module.
Form the induced $g(V)$-module $Ind(U)=U(g(V))\otimes_{U(g(V)_{-}+g(V)_{0})}U$.
Set $V(U)=Ind(U)/J(Ind(U))$. Then $V(U)$ is a lowest weight weak $V$-module.
If $V$ is rational, it follows from the complete reducibility of $V(U)$ that
$V(U)$ is irreducible.
The following is our generalized nuclear democracy theorem 
of Tsuchiya and Kanie [TK].

\bt{t4.14} 
Let $W_{1}$ and $W_{2}$ be $V$-modules.
Let $U$ be a $g(V)_{0}$-module on which $L(0)$ acts as a scalar 
$h$ and let $I_{0}(\cdot,x)$
be a linear injective map
{}from $U$ to $\left({\rm Hom}_{{\C}}(W_{1},W_{2})\right)\{x\}$
such that for any $u\in U$, $I_{0}(u,x)$ satisfies the truncation
condition $(G1)$, the $L(-1)$-bracket formula $(G2)$ and the
following condition: 
\begin{eqnarray}
& &(x_{1}-x_{2})^{{\rm wt}a-1}Y_{W_{2}}(a,x_{1})I_{0}(u,x_{2})-
(-x_{2}+x_{1})^{{\rm wt}a-1}I_{0}(u,x_{2})Y_{W_{1}}(a,x_{1})\nonumber\\
&=&x_{1}^{-1}\delta\left({x_{2}\over x_{1}}\right)I_{0}(a_{{\rm wt}a-1}u,x_{2})
\end{eqnarray}
for any $a\in V, u\in U$. Then there exists a lowest weight weak $V$-module 
$W$ with $U$ as its 
lowest weight subspace generating $W$ and there is a
unique intertwining operator $I(\cdot,x)$ of type 
$\left(\!\begin{array}{c}W_{2}\\ W W_{1}\end{array}\!\right)$
extending $I_{0}(\cdot,x)$. In particular, if $V$ is rational and $U$ is irreducible, $W$ is 
irreducible.
\et

{\bf Proof.} Since 
$\displaystyle{(x_{1}-x_{2})\delta\left({x_{2}\over x_{1}}\right)=0}$,
we have
\begin{eqnarray}
(x_{1}-x_{2})^{{\rm wt}a+i}Y_{W_{2}}(a,x_{1})I_{0}(u,x_{2})
=(-x_{2}+x_{1})^{{\rm wt}a+i}I_{0}(u,x_{2})Y_{W_{1}}(a,x_{1})
\end{eqnarray}
for $a\in V, u\in U,i\in {\N}$. Then by definition
$I_{0}(u,x)\in G(W_{1},W_{2})$ for any $u\in U$ and
\begin{eqnarray}
& &a_{m}\circ I_{0}(u,x)=0\;\;\;\mbox{for }m\ge {\rm wt }a,\\
& &a_{{\rm wt}a-1}\circ I_{0}(u,x)=I_{0}(a_{{\rm wt}a-1}u,x).
\end{eqnarray}
Set $\bar{U}:=\{I(u,x)|u\in U\}\subseteq G(W_{1},W_{2})$. Then
$\bar{U}$ is a $g(V)_{0}$-submodule of $G(W_{1},W_{2})$ and $\bar{U}$
as a $g(V)_{0}$-module is isomorphic to $U$.
Let $W=U(g(V))\bar{U}$ be the $V$ 
or $g(V)$-submodule of $G(W_{1},W_{2})$. Then $W=U(g_{+})\bar{U}$ is a
lower-truncated ${\Z}$-graded weak $V$-module generated by $U$.
Then we have a 
natural intertwining operator of type 
$\left(\!\begin{array}{c}W_{2}\\ W W_{1}\end{array}\!\right)$.
The uniqueness is clear. Then the proof is complete.$\;\;\;\;\Box$.

Let ${\bf g}$ be a finite-dimensional simple Lie algebra, let ${\bf h}$ be
a Cartan subalgebra, let $\Delta$ be the root system of
${\bf g}$ and let $\langle\cdot,\cdot\rangle$ be the normalized Killing
form on ${\bf g}$ [K]. For any
linear functional $\lambda\in {\bf h}^{*}$, we denote by $L(\lambda)$
the irreducible highest weight  ${\bf g}$-module with
highest weight $\lambda$. 

Let $\hat{{\bf g}}={\C}[t,t^{-1}]\otimes {\bf g}\oplus {\C}c$
be the corresponding affine Lie algebra and let $\tilde{{\bf
g}}=\hat{{\bf g}}\oplus {\C}d$ be the extended affine algebra.
For any $\ell\in {\C}, \lambda\in {\bf h}^{*}$, 
let $L(\ell,\lambda)$ be the irreducible highest  weight
$\hat{{\bf g}}$-module of level $\ell$.
For any ${\bf g}$-module $U$, let $\hat{U}$ be the loop $\hat{{\bf
g}}$-module ${\C}[t,t^{-1}]\otimes U$ of level $0$.
It is well known (cf. [FZ], [L1]) that each $L(\ell,0)$ has a
vertex operator algebra structure except when $\ell$ is the negative dual
Coxeter number. Then we have the following nuclear democracy theorem of 
Tsuchiya and Kanie.
(To be precise, this was proved only for ${\bf g}=sl_{2}$ in [TK].) 

\bp{p4.14}
Let $\ell$ be a positive integer and let 
$W_{2}, W_{3}$ be $L(\ell,0)$-modules.
Let $\lambda$ be a linear functional on ${\bf h}$, let $L(\lambda)$ be the 
irreducible highest weight ${\bf g}$-module with highest weight $\lambda$ and 
let $\Phi (\cdot,x)$ be a nonzero linear map
{}from $L(\lambda)$ to ${\rm Hom}(W_{2},W_{3})\{x\}$ such that
\begin{eqnarray}
& &[a(m), \Phi (u,x)]=x^{m}\Phi(a(0)u,x);\label{e4.38}\\
& &[L(-1), \Phi (u,x)]={d\over dx}\Phi (u,x)\label{e4.39}
\end{eqnarray} 
for any $a\in {\bf g}\subseteq L(\ell,0), u\in L(\lambda), m\in
{\Z}$,
Then $L(\ell,\lambda)$ is an irreducible $L(\ell,0)$-module 
and there is a unique intertwining operator
$I(\cdot,x)$ on $L(\ell,\lambda)$  in the sense of [FHL] extending $\Phi(\cdot,x)$.
\ep

{\bf Proof.} Writing (\ref{e4.38}) in terms of generating
functions, we obtain
\begin{eqnarray}
& &[Y(a,x_{1}),\Phi(u,x_{2})]=x_{2}^{-1}\delta\left({x_{1}\over x_{2}}\right)
\Phi(a_{0}u,x_{2}).
\end{eqnarray}
Since $\displaystyle{(x_{1}-x_{2})\delta\left({x_{1}\over
x_{2}}\right)=0}$, we get
\begin{eqnarray}\label{e4.44}
(x_{1}-x_{2})[Y(a,x_{1}),\Phi(u,x_{2})]=0
\end{eqnarray}
for any $a\in {\bf g}, u\in L(\lambda)$. Since ${\bf g}$
generates $L(\ell,0)$ as a vertex operator algebra, similar to
the proof of Proposition 4.5 it follows from the proof of Proposition 3.2.7 
in [L1] that $I(u,x)$ satisfies (G3) for
any $a\in L(\ell,0)$. Furthermore, (\ref{e4.38}) implies (G2).
Therefore $\Phi(u,x)\in G(W_{2},W_{3})$ for $u\in L(\lambda)$. From 
(\ref{e4.21}) and (\ref{e4.44}) we obtain
\begin{eqnarray}
& &a(0)\circ \Phi(u,x)=[a(0),\Phi(u,x)]=\Phi (a(0)u,x);\\
& &a(m)\circ \Phi (u,x)=0\;\;\;\;\mbox{for any }a\in {\bf g}, m>0,
u\in L(\lambda).
\end{eqnarray}
Then $\Phi$ is a ${\bf g}$-homomorphism. Consequently, $L(\lambda)$ is embedded
into $G(W_{2},W_{3})$ by $\Phi$. Let $W$ be the 
$V$-submodule generated by 
$L(\lambda)$. 
Then $W$ is a certain quotient module of $M(\ell,\lambda)$. From the rationality of
$L(\ell,0)$, we get
$W=L(\ell,\lambda)$. By Theorem \ref{t4.8}, we obtain an intertwining vertex
operator $I(\cdot,x)$ of type $\left(\!\begin{array}{c}W_{3})\\L(\ell,
\lambda)W_{2})\end{array}\!\right)$. The uniqueness is clear.
Then the proof is complete.$\;\;\;\;\Box$

\begin{rem}\label{r4.14}
Under the conditions of Proposition \ref{p4.14}, we obtain an intertwining operator $I(\cdot,x)$ 
of type $\left(\!\begin{array}{c}W_{3}\\L(\ell,\lambda) W_{2}\end{array}\!\right)$.
It follows from commutator formula (\ref{ec}) that
$$[L(m), I(u,x)]=x^{m}\left((m+1)h+x{d\over dx}\right)I(u,x)$$
for $u\in L(\lambda)$, where $h$ is the lowest weight of $L(\ell,\lambda)$. Thus
\begin{eqnarray}
[L(m), \Phi(u,x)]=x^{m}\left((m+1)h+x{d\over dx}\right)\Phi(u,x)\;\;\;\mbox{ for }u\in L(\lambda), m\in {\Z}.
\end{eqnarray}
\end{rem}

In many references the notion of loop $\hat{{\bf g}}$-module was used to 
define intertwining operators. Next we 
shall discuss this issue.

Suppose $L(\ell,\lambda_{i})$ $(i=1,2,3)$ are $L(\ell,0)$-modules. Let
$I(\cdot,x)$ be an intertwining operator of type
$\left(\!\begin{array}{c}L(\ell,\lambda_{3})\\L(\ell,\lambda_{1})
L(\ell,\lambda_{2})\end{array}\!\right)$. As before, we set
\begin{eqnarray}
I^{o}(u_{1},x)=x^{h_{1}+h_{2}-h_{3}}I(u_{1},x)=\sum_{n\in {\Z}}I_{u_{1}}(n)x^{-n-1}
\;\;\;\mbox{ for any }u_{1}\in L(\ell,\lambda_{1}).
\end{eqnarray}
Then (the second identity follows from Proposition 2.5)
\begin{eqnarray}
& &[a(m), I_{u}(n)]=I_{au}(m+n);\label{e4.40}\\
& &[L(0), I_{u}(n)]=(h_{3}-h_{2}-n-1)I_{u}(n)\label{e4.41}
\end{eqnarray}
for $a\in {\bf g}, u\in L(\lambda_{1}), m,n\in {\Z}.$
Then $I^{o}(\cdot,x)$ naturally gives rise to a linear map $R_{I}$ from 
${\C}[t,t^{-1}]\otimes L(\lambda_{1})\otimes L(\ell,\lambda_{2})$ to
$L(\ell,\lambda_{3})$ such that
$$R_{I}(t^{n}\otimes u_{1}\otimes u_{2})=I_{u_{1}}(n)u_{2}\;\;\;\mbox{for }
u_{1}\in L(\lambda_{1}),u_{2}\in L(\ell,\lambda_{2}),n\in {\Z}.$$
(\ref{e4.40}) is equivalent to say that
the map $R_{I}$ is a  $\hat{{\bf g}}$-homomorphism from
$\hat{L}(\lambda_{1})\otimes L(\ell,\lambda_{2})$ to
$L(\ell,\lambda_{3})$. From (\ref{e4.41}) we get
\begin{eqnarray}
L(0)(t^{n}\otimes u_{1}\otimes u_{2})=t^{n}\otimes u_{1}\otimes L(0)u_{2}
+(h_{3}-h_{2}-n-1)(t^{n}\otimes u_{1}\otimes u_{2})
\end{eqnarray}
for $u_{1}\in L(\lambda_{1}),u_{2}\in L(\ell,\lambda_{2}), n\in {\Z}$. Then
\begin{eqnarray}\label{eda}
(h_{3}-L(0))(t^{n}\otimes u_{1}\otimes u_{2})=t^{n}\otimes u_{1}\otimes (
h_{2}-L(0))u_{2}+(n+1)(t^{n}\otimes u_{1}\otimes u_{2}).
\end{eqnarray}

View $\hat{L}(\lambda_{1})$ as a $\tilde{{\bf g}}$-module with
$d=(1+t{d\over dt})\otimes 1$ and view $L(\ell,\lambda)$ as a
$\tilde{{\bf g}}$-module with $d=h-L(0)$ where $h$ is the lowest weight.
Then it follows from (\ref{eda}) that $R_{I}$ is a $\tilde{{\bf
g}}$-homomorphism. Then we obtain a linear map:
\begin{eqnarray}
R :& &I\!\left(\!\begin{array}{c}L(\ell,\lambda_{3})\\L(\ell,\lambda_{1})
L(\ell,\lambda_{2})\end{array}\!\right)\rightarrow {\rm
Hom}_{\hat{{\bf g}}}(\hat{L}(\lambda_{1})\otimes
L(\ell,\lambda_{2}), L(\ell,\lambda_{3}));\nonumber\\
& &I(\cdot,x)\mapsto R_{I}.
\end{eqnarray}

In some references, an intertwining operator of type
$\left(\!\begin{array}{c}L(\ell,\lambda_{3})\\L(\ell,\lambda_{1})
L(\ell,\lambda_{2})\end{array}\!\right)$ is defined to be a 
$\tilde{{\bf g}}$-module homomorphism 
{}from $\hat{L}(\lambda_{1})\otimes L(\ell,\lambda_{2})$ to 
$L(\ell,\lambda_{3})$. The following
proposition asserts that this definition is equivalent to FHL's definition.

\bp{p4.15}
The intertwining operator space 
$I\!\left(\!\begin{array}{c}L(\ell,\lambda_{3})\\
L(\ell,\lambda_{1})
L(\ell,\lambda_{2})\end{array}\!\right)$ is naturally isomorphic to 
the space of $\tilde{{\bf g}}$-homomorphisms from 
$\hat{L}(\lambda)\otimes L(\ell,\lambda_{2})$ 
to $L(\ell,\lambda_{3})$.
\ep

{\bf Proof.} From the above discussion we see that for any intertwining 
operator $I(\cdot,x)$
of type $\left(\!\begin{array}{c}L(\ell,\lambda_{3})\\L(\ell,\lambda_{1})
L(\ell,\lambda_{2})\end{array}\!\right)$, we obtain a 
$\tilde{{\bf g}}$-homomorphism 
$R_{I}$.
Conversely, let $f$ be a $\tilde{{\bf g}}$-homomorphism from $\hat{L}(\lambda)
\otimes L(\ell,\lambda_{2})$ 
to $L(\ell,\lambda_{3})$. Then we define a linear map $\Phi(\cdot,x)$ from 
$L(\lambda_{1})$ to 
${\rm Hom}(L(\ell,\lambda_{2}),L(\ell,\lambda_{2}))\{x\}$ such that
\begin{eqnarray}
\Phi(u_{1},x)u_{2}=x^{h_{3}-h_{1}-h_{2}}\sum_{n\in {\Z}}f(t^{n}\otimes u_{1}
\otimes u_{2})x^{-n-1}
\end{eqnarray}
for $u_{1}\in L(\lambda_{1}),u_{2}\in L(\ell,\lambda_{2})$. Then 
$\Phi(\cdot,x)$ satisfies (\ref{e4.38})
and
\begin{eqnarray}\label{e4.45}
[L(0), \Phi(u_{1},x)]=\left(h_{1}+x{d\over dx}\right)\Phi(u_{1},x)\;\;\;
\mbox{ for }u_{1}\in 
L(\lambda_{1}).
\end{eqnarray}
Then $\Phi(u_{1},x)\in E(L(\ell,\lambda_{2}), L(\ell,\lambda_{3}))$ for any 
$u_{1}\in L(\lambda_{1})$.
Similar to the proof of Proposition \ref{p4.14}, $L(\ell,\lambda_{1})$ is a 
submodule of
$E(L(\ell,\lambda_{2}), L(\ell,\lambda_{3}))$ generated by $L(\lambda_{1})$ and
there is a weak intertwining operator
$I(\cdot,x)$ from $L(\ell,\lambda)$ to ${\rm Hom}(L(\ell,\lambda_{2}),
L(\ell,\lambda_{3}))\{x\}$.
It is well known (cf. [HL0-4], [FLM]) that under the commutator formula (\ref{ec}), 
the $L(-1)$-bracket formula (I2) is equivalent to 
the $L(0)$-bracket formula$.\!$ Thus $\Phi(u_{1},x)\!\in \!
G(L(\ell,\lambda_{2}),L(\ell,\lambda_{3}))$
for $u_{1}\in L(\lambda_{1})$. Since $L(\lambda_{1})$ generates 
$L(\ell,\lambda_{1})$ by 
$U(\hat{{\bf g}})$, it follows from Proposition \ref{p4.5} that
$L(\ell,\lambda_{1})\subseteq G(L(\ell,\lambda_{2}),L(\ell,\lambda_{3}))$. Thus
$I(\cdot,x)$ is an intertwining operator. Then the proof is complete.$\;\;\;\;\Box$

Let $L(c,h)$ be the irreducible module of the Virasoro algebra Vir with 
central charge $c$ and lowest weight $h$. It is 
well known (cf. [FZ], [H1], [L1]) that $L(c,0)$ is a vertex operator algebra. 
Suppose that $L(c,h_{1})$ and $L(c,h_{2})$ are two
modules for the vertex operator algebra $L(c,0)$. Let $\Phi(x)\in
\left({\rm Hom}_{{\C}}(L(c,h_{1}),L(c,h_{2}))\right)\{x\}$ such that
\begin{eqnarray}\label{e4.50}
[L(m),\Phi(x)]=x^{m}\left((m+1)h+x{d\over dx}\right)\Phi(x)
\end{eqnarray}
for some complex number $h$. That is,
\begin{eqnarray}\label{e4.51}
[Y(\omega,x_{1}),\Phi(x_{2})]=x_{1}^{-1}\delta\left({x_{2}\over x_{1}}\right)
{d\over dx_{2}}\Phi(x_{2})
+hx_{1}^{-2}\delta'\left({x_{2}\over x_{1}}\right)\Phi(x_{2}).
\end{eqnarray}
Similarly to the proof of Proposition \ref{p4.12} we get 
$\Phi(x)\in G(L(c,h_{1}),L(c,h_{2}))$ and $\Phi(x)$ generates a $L(c,0)$-module
$M$ which is a lowest weight Virasoro algebra module with lowest weight $h$ in ${\rm
G}(L(c,h_{1}),L(c,h_{2}))$. If $c=1-\frac{6(p-q)^{2}}{pq}$, where 
$p,q\in \{2,3\cdots\}$ are relatively prime,
$L(c,0)$ is rational ([DMZ], [W]). Therefore
$M=L(c,h)$. Then we obtain an intertwining vertex operator of type
$\left(\!\begin{array}{c}L(c,h_{2})\\L(c,h) L(c,h_{2})\end{array}\!\right)$. Thus 
we have

\bp{pmin}
If $c=1-\frac{6(p-q)^{2}}{pq}$, where 
$p,q\in \{2,3\cdots\}$ are relatively prime, let $L(c,h_{1}), L(c,h_{2})$ be 
$L(c,0)$-modules
and let $\Phi(x)$ satisfy (\ref{e4.50}). Then there exists a unique
intertwining vertex operator of type
$\left(\!\begin{array}{c}L(c,h_{2})\\L(c,h) L(c,h_{2})\end{array}\!\right)$
extending $\Phi(x)$.
\ep

\section{Appendix}
The main purpose of this appendix is to give an example to show that
the generalized form of the nuclear democracy 
theorem may not be true if $V$ 
is not rational. We use the same notions as in Section 4.
 Let $\ell$ be a positive integer and let 
${\C}_{\ell}$ be the $({\C}[t]\otimes {\bf g}+{\C}c)$-module such that
$c$ acts as $\ell$ and ${\C}[t]\otimes {\bf g}$ acts as zero. Set
$$M(\ell,{\C})=U({\bf g})_{U({\C}[t]\oplus {\bf g}+{\C}c)}{\C}_{\ell}.$$
Then $M(\ell,{\C})$
 is a vertex operator algebra and any restricted $\hat{{\bf g}}$-modules of 
level $\ell$ is a 
$M(\ell,{\C})$-module (cf. [FZ], [L1]). Consequently, $M(\ell,{\C})$ is 
irrational.
Since $L(\ell,0)$ is rational, we may choose an $\lambda$ such that 
$L(\ell,\lambda)$ is not a $L(\ell,0)$-module.
Let $\Phi (x)$ be the identity map from $L(\ell,\lambda)$ to $L(\ell,\lambda)$.
Then $\Phi$ satisfies all the conditions in Theorem \ref{t4.14}.
If  we could extend $\Phi$ to an intertwining operator on $L(\ell,0)$, then
we would have an intertwining operator of type
$\left(\!\begin{array}{c}L(\ell,\lambda)\\L(\ell,0)L(\ell,\lambda)\end{array}\!
\right)$ so that $L(\ell,\lambda)$ would be a $L(\ell,0)$-module.
This would contradict the assumption that $L(\ell,\lambda)$ is not a 
$L(\ell,0)$-module.

\end{document}